\documentclass{amsart}
\usepackage{amssymb}

\usepackage{amscd}
\usepackage{thmdefs}

\input{tcilatex}
\theoremstyle{definition}
\theoremstyle{remark}
\numberwithin{equation}{section}

\input tcilatex

\begin{document}
\title[Toeplitz Matricial R-transform Theory]{Amalgamated R-transform Theory on Toeplitz-Probability Spaces over Toeplitz
Matricial Algebras }
\author{Ilwoo Cho}
\address{Dep. of Math, Univ. of Iowa, Iowa City, IA, U.S.A}
\email{icho@math.uiowa.edu}
\maketitle

\begin{abstract}
In this paper, we will consider a noncommutative probability space, $(T,E),$
over a Toeplitz matricial algebra $B=\mathcal{C}^{N},$ for $N\in \Bbb{N}$,
induced by a (scalar-valued) noncommutative probability space, $(A,\mathcal{%
\varphi }),$ with a suitable $B$-functional, $E:T\rightarrow B,$ defined by $%
\varphi .$ On this framework, we will observe amalgamated R-transform
calculus with respect to a $B$-functional, $E.$ The technique and idea are
came from those in [12] and Free Probability of type $B$ studied in [31].
\end{abstract}

\strut

\strut

\strut

Voiculescu developed Free Probability Theory. Here, the classical concept of
Independence in Probability theory is replaced by a noncommutative analogue
called Freeness (See [9]). There are two approaches to study Free
Probability Theory. One of them is the original analytic approach of
Voiculescu and the other one is the combinatorial approach of Speicher and
Nica (See [1], [2] and [3]).\medskip

To observe the free additive convolution and free multiplicative convolution
of two distributions induced by free random variables in a noncommutative
probability space (over $B=\Bbb{C}$), Voiculescu defined the R-transform and
the S-transform, respectively. These show that to study distributions is to
study certain ($B$-)formal series for arbitrary noncommutative
indeterminants.\strut

\strut \strut

Speicher defined the free cumulants which are the main objects in the
combinatorial approach of Free Probability Theory. And he developed free
probability theory by using Combinatorics and Lattice theory on collections
of noncrossing partitions (See [3]). Also, Speicher considered the
operator-valued free probability theory, which is also defined and observed
analytically by Voiculescu, when $\Bbb{C}$ is replaced to an arbitrary
algebra $B$ (See [1]). Nica defined R-transforms of several random variables
(See [2]). He defined these R-transforms as multivariable formal series in
noncommutative several indeterminants. To observe the R-transform, the
M\"{o}bius Inversion under the embedding of lattices plays a key role (See
[1],[3],[5],[12],[13] and [17]).\strut

\strut

In [12], [19] and [20], we observed the amalgamated R-transform calculus.
Actually, amalgamated R-transforms are defined originally by Voiculescu (See
[10]) and are characterized combinatorially by Speicher (See [1]). In [12],
we defined amalgamated R-transforms slightly differently from those defined
in [1] and [10]. We defined them as $B$-formal series and tried to
characterize, like in [2] and [3]. The main tool which is considered, for
studying amalgamated R-transforms is amalgamated boxed convolution. After
defining boxed convolution over an arbitrary algebra $B,$ we could get that

\strut

\begin{center}
$R_{x_{1},...,x_{s}}\,\,\frame{*}\,%
\,R_{y_{1},...,y_{s}}^{symm(1_{B})}=R_{x_{1}y_{1},...,x_{s}y_{s}},$ for any $%
s\in \Bbb{N},$
\end{center}

\strut

where $x_{j}$ and $y_{j}$ are free $B$-valued random variables.

\strut 

In this paper, we will consider a noncommutative probability space, $(T,E),$
over a Toeplitz matricial algebra $B=\mathcal{C}^{N},$ for $N\in \Bbb{N}$,
induced by a (scalar-valued) noncommutative probability space, $(A,\mathcal{%
\varphi }),$ with a suitable $B$-functional, $E:T\rightarrow B,$ defined by $%
\varphi .$ On this framework, we will observe amalgamated R-transform
calculus with respect to a $B$-functional, $E.$ The technique and idea are
came from those in [12] and Free Probability of type $B$ studied in [31]. By
the construction of this amalgamated noncommutative probability space, we
can easily recognize that to compute the partition-dependent moments of
random variables, we do not need to concern about the insertion property. We
will concentrate on observing the moment series and R-transforms of random
variables in this structure and after considering the (combinatorial)\strut 
amalgamated freeness, we construct the R-transform calculus.

\strut

\strut

\strut

\section{Preliminaries}

\strut

\strut

\subsection{Free Probability Theory}

\strut

\strut

\strut \strut

In this section, we will summarize and introduced the basic results from [1]
and [12]. Throughout this section, let $B$ be a unital algebra. The
algebraic pair $(A,\varphi )$ is said to be a noncommutative probability
space over $B$ (shortly, NCPSpace over $B$) if $A$ is an algebra over $B$
(i.e $1_{B}=1_{A}\in B\subset A$) and $\varphi :A\rightarrow B$ is a $B$%
-functional (or a conditional expectation) ; $\varphi $ satisfies

\strut

\begin{center}
$\varphi (b)=b,$ for all $b\in B$
\end{center}

and

\begin{center}
$\varphi (bxb^{\prime })=b\varphi (x)b^{\prime },$ for all $b,b^{\prime }\in
B$ and $x\in A.$
\end{center}

\strut

\strut Let $(A,\varphi )$ be a NCPSpace over $B.$ Then, for the given $B$%
-functional, we can determine a moment multiplicative function $\widehat{%
\varphi }=(\varphi ^{(n)})_{n=1}^{\infty }\in I(A,B),$ where

\strut

\begin{center}
$\varphi ^{(n)}(a_{1}\otimes ...\otimes a_{n})=\varphi (a_{1}....a_{n}),$
\end{center}

\strut \strut

for all $a_{1}\otimes ...\otimes a_{n}\in A^{\otimes _{B}n},$ $\forall n\in 
\Bbb{N}.$

\strut

\strut We will denote noncrossing partitions over $\{1,...,n\}$ ($n\in \Bbb{N%
}$) by $NC(n).$ Define an ordering on $NC(n)$ ;

\strut

$\theta =\{V_{1},...,V_{k}\}\leq \pi =\{W_{1},...,W_{l}\}$\strut $\overset{%
def}{\Leftrightarrow }$ For each block $V_{j}\in \theta $, there exists only
one block $W_{p}\in \pi $ such that $V_{j}\subset W_{p},$ for $j=1,...,k$
and $p=1,...,l.$

\strut

Then $(NC(n),\leq )$ is a complete lattice with its minimal element $%
0_{n}=\{(1),...,(n)\}$ and its maximal element $1_{n}=\{(1,...,n)\}$. We
define the incidence algebra $I_{2}$ by a set of all complex-valued\
functions $\eta $ on $\cup _{n=1}^{\infty }\left( NC(n)\times NC(n)\right) $
satisfying $\eta (\theta ,\pi )=0,$ whenever $\theta \nleq \pi .$ Then,
under the convolution

\strut

\begin{center}
$*:I_{2}\times I_{2}\rightarrow \Bbb{C}$
\end{center}

defined by

\begin{center}
$\eta _{1}*\eta _{2}(\theta ,\pi )=\underset{\theta \leq \sigma \leq \pi }{%
\sum }\eta _{1}(\theta ,\sigma )\cdot \eta _{2}(\sigma ,\pi ),$
\end{center}

\strut

$I_{2}$ is indeed an algebra of complex-valued functions. Denote zeta,
M\"{o}bius and delta functions in the incidence algebra $I_{2}$ by $\zeta ,$ 
$\mu $ and $\delta ,$ respectively. i.e

\strut

\begin{center}
$\zeta (\theta ,\pi )=\left\{ 
\begin{array}{lll}
1 &  & \theta \leq \pi \\ 
0 &  & otherwise,
\end{array}
\right. $
\end{center}

\strut

\begin{center}
$\delta (\theta ,\pi )=\left\{ 
\begin{array}{lll}
1 &  & \theta =\pi \\ 
0 &  & otherwise,
\end{array}
\right. $
\end{center}

\strut

and $\mu $ is the ($*$)-inverse of $\zeta .$ Notice that $\delta $ is the ($%
* $)-identity of $I_{2}.$ By using the same notation ($*$), we can define a
convolution between $I(A,B)$ and $I_{2}$ by

\strut

\begin{center}
$\widehat{f}\,*\,\eta \left( a_{1},...,a_{n}\,;\,\pi \right) =\underset{\pi
\in NC(n)}{\sum }\widehat{f}(\pi )(a_{1}\otimes ...\otimes a_{n})\eta (\pi
,1_{n}),$
\end{center}

\strut

where $\widehat{f}\in I(A,B)$, $\eta \in I_{1},$ $\pi \in NC(n)$ and $%
a_{j}\in A$ ($j=1,...,n$), for all $n\in \Bbb{N}.$ Notice that $\widehat{f}%
*\eta \in I(A,B),$ too. Let $\widehat{\varphi }$ be a moment multiplicative
function in $I(A,B)$ which we determined before. Then we can naturally
define a cumulant multiplicative function $\widehat{c}=(c^{(n)})_{n=1}^{%
\infty }\in I(A,B)$ by

\strut

\begin{center}
$\widehat{c}=\widehat{\varphi }*\mu $ \ \ \ or \ \ $\widehat{\varphi }=%
\widehat{c}*\zeta .$
\end{center}

\strut

\strut

This says that if we have a moment\ multiplicative function, then we always
get a cumulant multiplicative function and vice versa, by $(*).$ This
relation is so-called ''M\"{o}bius Inversion''. More precisely, we have

\strut

\begin{center}
$
\begin{array}{ll}
\varphi (a_{1}...a_{n}) & =\varphi ^{(n)}(a_{1}\otimes ...\otimes a_{n}) \\ 
& =\underset{\pi \in NC(n)}{\sum }\widehat{c}(\pi )(a_{1}\otimes ...\otimes
a_{n})\zeta (\pi ,1_{n}) \\ 
& =\underset{\pi \in NC(n)}{\sum }\widehat{c}(\pi )(a_{1}\otimes ...\otimes
a_{n}),
\end{array}
$
\end{center}

\strut

for all $a_{j}\in A$ and $n\in \Bbb{N}.$ Or equivalently,

\strut

\begin{center}
$
\begin{array}{ll}
c^{(n)}(a_{1}\otimes ...\otimes a_{n}) & =\underset{\pi \in NC(n)}{\sum }%
\widehat{\varphi }(\pi )(a_{1}\otimes ...\otimes a_{n})\mu (\pi ,1_{n}).
\end{array}
$
\end{center}

\strut \strut

Now, let $(A_{i},\varphi _{i})$ be NCPSpaces over $B,$ for all $i\in I.$
Then we can define a amalgamated free product of $A_{i}$ 's and amalgamated
free product of $\varphi _{i}$'s by

\begin{center}
$A\equiv *_{B}A_{i}$ \ \ and \ $\varphi \equiv *_{i}\varphi _{i},$
\end{center}

respectively. Then, by Voiculescu, $(A,\varphi )$ is again a NCPSpace over $%
B $ and, as a vector space, $A$ can be represented by

\begin{center}
\strut $A=B\oplus \left( \oplus _{n=1}^{\infty }\left( \underset{i_{1}\neq
...\neq i_{n}}{\oplus }(A_{i_{1}}\ominus B)\otimes ...\otimes
(A_{i_{n}}\ominus B)\right) \right) ,$
\end{center}

where $A_{i_{j}}\ominus B=\ker \varphi _{i_{j}}.$ We will use Speicher's
combinatorial definition of amalgamated free product of $B$-functionals ;

\strut

\begin{definition}
Let $(A_{i},\varphi _{i})$ be NCPSpaces over $B,$ for all $i\in I.$ Then $%
\varphi =*_{i}\varphi _{i}$ is the amalgamated free product of $B$%
-functionals $\varphi _{i}$'s on $A=*_{B}A_{i}$ if the cumulant
multiplicative function $\widehat{c}=\widehat{\varphi }*\mu \in I(A,B)$ has
its restriction to $\underset{i\in I}{\cup }A_{i},$ $\underset{i\in I}{%
\oplus }\widehat{c_{i}},$ where $\widehat{c_{i}}$ is the cumulant
multiplicative function induced by $\varphi _{i},$ for all $i\in I.$ i.e

\strut

\begin{center}
$c^{(n)}(a_{1}\otimes ...\otimes a_{n})=\left\{ 
\begin{array}{lll}
c_{i}^{(n)}(a_{1}\otimes ...\otimes a_{n}) &  & \text{if }\forall a_{j}\in
A_{i} \\ 
0_{B} &  & otherwise.
\end{array}
\right. $
\end{center}
\end{definition}

\strut

Now, we will observe the freeness over $B$ ;

\strut

\begin{definition}
Let $(A,\varphi )$\strut be a NCPSpace over $B.$

\strut

(1) Subalgebras containing $B,$ $A_{i}\subset A$ ($i\in I$) are free (over $%
B $) if we let $\varphi _{i}=\varphi \mid _{A_{i}},$ for all $i\in I,$ then $%
*_{i}\varphi _{i}$ has its cumulant multiplicative function $\widehat{c}$
such that its restriction to $\underset{i\in I}{\cup }A_{i}$ is $\underset{%
i\in I}{\oplus }\widehat{c_{i}},$ where $\widehat{c_{i}}$ is the cumulant
multiplicative function induced by each $\varphi _{i},$ for all $i\in I.$

\strut

(2) Sebsets $X_{i}$ ($i\in I$) are free (over $B$) if subalgebras $A_{i}$'s
generated by $B$ and $X_{i}$'s are free in the sense of (1). i.e If we let $%
A_{i}=A\lg \left( X_{i},B\right) ,$ for all $i\in I,$ then $A_{i}$'s are
free over $B.$
\end{definition}

\strut

In [1], Speicher showed that the above combinatorial freeness with
amalgamation can be used alternatively with respect to Voiculescu's original
freeness with amalgamation.

\strut

Let $(A,\varphi )$ be a NCPSpace over $B$ and let $x_{1},...,x_{s}$ be $B$%
-valued random variables ($s\in \Bbb{N}$). Define $(i_{1},...,i_{n})$-th
moment of $x_{1},...,x_{s}$ by

\strut

\begin{center}
$\varphi (x_{i_{1}}b_{i_{2}}x_{i_{2}}...b_{i_{n}}x_{i_{n}}),$
\end{center}

\strut

for arbitrary $b_{i_{2}},...,b_{i_{n}}\in B,$ where $(i_{1},...,i_{n})\in
\{1,...,s\}^{n},$ $\forall n\in \Bbb{N}.$ Similarly, define a symmetric $%
(i_{1},...,i_{n})$-th moment by the fixed $b_{0}\in B$ by

\strut

\begin{center}
$\varphi (x_{i_{1}}b_{0}x_{i_{2}}...b_{0}x_{i_{n}}).$
\end{center}

\strut

If $b_{0}=1_{B},$ then we call this symmetric moments, trivial moments.

\strut

Cumulants defined below are main tool of combinatorial free probability
theory ; in [12], we defined the $(i_{1},...,i_{n})$-th cumulant of $%
x_{1},...,x_{s}$ by

\strut

\begin{center}
$k_{n}(x_{i_{1}},...,x_{i_{n}})=c^{(n)}(x_{i_{1}}\otimes
b_{i_{2}}x_{i_{2}}\otimes ...\otimes b_{i_{n}}x_{i_{n}}),$
\end{center}

\strut

for $b_{i_{2}},...,b_{i_{n}}\in B,$ arbitrary, and $(i_{1},...,i_{n})\in
\{1,...,s\}^{n},$ $\forall n\in \Bbb{N},$ where $\widehat{c}%
=(c^{(n)})_{n=1}^{\infty }$ is the cumulant multiplicative function induced
by $\varphi .$ Notice that, by M\"{o}bius inversion, we can always take such 
$B$-value whenever we have $(i_{1},...,i_{n})$-th moment of $%
x_{1},...,x_{s}. $ And, vice versa, if we have cumulants, then we can always
take moments. Hence we can define a symmetric $(i_{1},...,i_{n})$-th\
cumulant by $b_{0}\in B$ of $x_{1},...,x_{s}$ by

\strut

\begin{center}
$k_{n}^{symm(b_{0})}(x_{i_{1}},...,x_{i_{n}})=c^{(n)}(x_{i_{1}}\otimes
b_{0}x_{i_{2}}\otimes ...\otimes b_{0}x_{i_{n}}).$
\end{center}

\strut

If $b_{0}=1_{B},$ then it is said to be trivial cumulants of $%
x_{1},...,x_{s} $.

\strut

By Speicher, it is shown that subalgebras $A_{i}$ ($i\in I$) are free over $%
B $ if and only if all mixed cumulants vanish.

\strut

\begin{proposition}
(See [1] and [12]) Let $(A,\varphi )$ be a NCPSpace over $B$ and let $%
x_{1},...,x_{s}\in (A,\varphi )$ be $B$-valued random variables ($s\in \Bbb{N%
}$). Then $x_{1},...,x_{s}$ are free if and only if all their mixed
cumulants vanish. $\square $
\end{proposition}

\strut

\strut

\strut

\strut

\strut

\subsection{Amalgamated R-transform Theory}

\strut

\strut

In this section, we will define an R-transform of several $B$-valued random
variables. Note that to study R-transforms is to study operator-valued
distributions. R-transforms with single variable is defined by Voiculescu
(over $B,$ in particular, $B=\Bbb{C}$. See [9] and [10]). Over $\Bbb{C},$
Nica defined multi-variable R-transforms in [2]. In [12], we extended his
concepts, over $B.$ R-transforms of $B$-valued random variables can be
defined as $B$-formal series with its $(i_{1},...,i_{n})$-th coefficients, $%
(i_{1},...,i_{n})$-th cumulants of $B$-valued random variables, where $%
(i_{1},...,i_{n})\in \{1,...,s\}^{n},$ $\forall n\in \Bbb{N}.$

\strut

\begin{definition}
Let $(A,\varphi )$ be a NCPSpace over $B$ and let $x_{1},...,x_{s}\in
(A,\varphi )$ be $B$-valued random variables ($s\in \Bbb{N}$). Let $%
z_{1},...,z_{s}$ be noncommutative indeterminants. Define a moment series of 
$x_{1},...,x_{s}$, as a $B$-formal series, by

\strut

\begin{center}
$M_{x_{1},...,x_{s}}(z_{1},...,z_{s})=\sum_{n=1}^{\infty }\underset{%
i_{1},..,i_{n}\in \{1,...,s\}}{\sum }\varphi
(x_{i_{1}}b_{i_{2}}x_{i_{2}}...b_{i_{n}}x_{i_{n}})\,z_{i_{1}}...z_{i_{n}},$
\end{center}

\strut

where $b_{i_{2}},...,b_{i_{n}}\in B$ are arbitrary for all $%
(i_{2},...,i_{n})\in \{1,...,s\}^{n-1},$ $\forall n\in \Bbb{N}.$

\strut

Define an R-transform of $x_{1},...,x_{s}$, as a $B$-formal series, by

\strut

\begin{center}
$R_{x_{1},...,x_{s}}(z_{1},...,z_{s})=\sum_{n=1}^{\infty }\underset{%
i_{1},...,i_{n}\in \{1,...,s\}}{\sum }k_{n}(x_{i_{1}},...,x_{i_{n}})%
\,z_{i_{1}}...z_{i_{n}},$
\end{center}

\strut with

\begin{center}
$k_{n}(x_{i_{1}},...,x_{i_{n}})=c^{(n)}(x_{i_{1}}\otimes
b_{i_{2}}x_{i_{2}}\otimes ...\otimes b_{i_{n}}x_{i_{n}}),$
\end{center}

\strut

where $b_{i_{2}},...,b_{i_{n}}\in B$ are arbitrary for all $%
(i_{2},...,i_{n})\in \{1,...,s\}^{n-1},$ $\forall n\in \Bbb{N}.$ Here, $%
\widehat{c}=(c^{(n)})_{n=1}^{\infty }$ is a cumulant multiplicative function
induced by $\varphi $ in $I(A,B).$
\end{definition}

\strut

Denote a set of all $B$-formal series with $s$-noncommutative indeterminants
($s\in \Bbb{N}$), by $\Theta _{B}^{s}$. i.e if $g\in \Theta _{B}^{s},$ then

\begin{center}
$g(z_{1},...,z_{s})=\sum_{n=1}^{\infty }\underset{i_{1},...,i_{n}\in
\{1,...,s\}}{\sum }b_{i_{1},...,i_{n}}\,z_{i_{1}}...z_{i_{n}},$
\end{center}

where $b_{i_{1},...,i_{n}}\in B,$ for all $(i_{1},...,i_{n})\in
\{1,...,s\}^{n},$ $\forall n\in \Bbb{N}.$ Trivially, by definition, $%
M_{x_{1},...,x_{s}},$ $R_{x_{1},...,x_{s}}\in \Theta _{B}^{s}.$ By $\mathcal{%
R}_{B}^{s},$\ we denote a set of all R-transforms of $s$-$B$-valued random
variables. Recall that, set-theoratically,

\begin{center}
$\Theta _{B}^{s}=\mathcal{R}_{B}^{s},$ sor all $s\in \Bbb{N}.$
\end{center}

By the previous section, we can also define symmetric moment series and
symmetric R-transform by $b_{0}\in B,$ by

\strut

\begin{center}
$M_{x_{1},...,x_{s}}^{symm(b_{0})}(z_{1},...,z_{s})=\sum_{n=1}^{\infty }%
\underset{i_{1},...,i_{n}\in \{1,...,s\}}{\sum }\varphi
(x_{i_{1}}b_{0}x_{i_{2}}...b_{0}x_{i_{n}})\,z_{i_{1}}...z_{i_{n}}$
\end{center}

and

\begin{center}
$R_{x_{1},...,x_{s}}^{symm(b_{0})}(z_{1},...,z_{s})=\sum_{n=1}^{\infty }%
\underset{i_{1},..,i_{n}\in \{1,...,s\}}{\sum }%
k_{n}^{symm(b_{0})}(x_{i_{1}},...,x_{i_{n}})\,z_{i_{1}}...z_{i_{n}},$
\end{center}

with

\begin{center}
$k_{n}^{symm(b_{0})}(x_{i_{1}},...,x_{i_{n}})=c^{(n)}(x_{i_{1}}\otimes
b_{0}x_{i_{2}}\otimes ...\otimes b_{0}x_{i_{n}}),$
\end{center}

for all $(i_{1},...,i_{n})\in \{1,...,s\}^{n},$ $\forall n\in \Bbb{N}.$

\strut

If $b_{0}=1_{B},$ then we have trivial moment series and trivial R-transform
of $x_{1},...,x_{s}$ denoted by $M_{x_{1},...,x_{s}}^{t}$ and $%
R_{x_{1},...,x_{s}}^{t},$ respectively.

\strut

The followings are known in [1] and [12] ;

\strut

\begin{proposition}
Let $(A,\varphi )$ be a NCPSpace over $B$ and let $%
x_{1},...,x_{s},y_{1},...,y_{p}\in (A,\varphi )$ be $B$-valued random
variables, where $s,p\in \Bbb{N}.$ Suppose that $\{x_{1},...,x_{s}\}$ and $%
\{y_{1},...,y_{p}\}$ are free in $(A,\varphi ).$ Then

\strut

(1) $
R_{x_{1},...,x_{s},y_{1},...,y_{p}}(z_{1},...,z_{s+p})=R_{x_{1},...,x_{s}}(z_{1},...,z_{s})+R_{y_{1},...,y_{p}}(z_{1},...,z_{p}). 
$

\strut

(2) If $s=p,$ then $R_{x_{1}+y_{1},...,x_{s}+y_{s}}(z_{1},...,z_{s})=\left(
R_{x_{1},...,x_{s}}+R_{y_{1},...,y_{s}}\right) (z_{1},...,z_{s}).$

$\square $
\end{proposition}

\strut

The above proposition is proved by the characterization of freeness with
respect to cumulants. i.e $\{x_{1},...,x_{s}\}$ and $\{y_{1},...,y_{p}\}$
are free in $(A,\varphi )$ if and only if their mixed cumulants vanish. Thus
we have

\strut

$k_{n}(p_{i_{1}},...,p_{i_{n}})=c^{(n)}(p_{i_{1}}\otimes
b_{i_{2}}p_{i_{2}}\otimes ...\otimes b_{i_{n}}p_{i_{n}})$

$\ \ \ \ =\left( \widehat{c_{x}}\oplus \widehat{c_{y}}\right)
^{(n)}(p_{i_{1}}\otimes b_{i_{2}}p_{i_{2}}\otimes ...\otimes
b_{i_{n}}p_{i_{n}})$

$\ \ \ \ =\left\{ 
\begin{array}{lll}
k_{n}(x_{i_{1}},...,x_{i_{n}}) &  & or \\ 
k_{n}(y_{i_{1}},...,y_{i_{n}}) &  & 
\end{array}
\right. $

\strut

and if $s=p,$ then

$k_{n}(x_{i_{1}}+y_{i_{1}},...,x_{i_{n}}+y_{i_{n}})=c^{(n)}\left(
(x_{i_{1}}+y_{i_{1}})\otimes b_{i_{2}}(x_{i_{2}}+y_{i_{2}})\otimes
...\otimes b_{i_{n}}(x_{i_{n}}+y_{i_{n}})\right) $

$\ =c^{(n)}(x_{i_{1}}\otimes b_{i_{2}}x_{i_{2}}\otimes ...\otimes
b_{i_{n}}x_{i_{n}})+c^{(n)}(y_{i_{1}}\otimes b_{i_{2}}y_{i_{2}}\otimes
...\otimes b_{i_{n}}y_{i_{n}})+[Mixed]$

\strut

where $[Mixed]$ is the sum of mixed cumulants, by the bimodule map property
of $c^{(n)}$

\strut

$\ =k_{n}(x_{i_{1}},...,x_{i_{n}})+k_{n}(y_{i_{1}},...,y_{i_{n}})+0_{B}.$

\strut

Now, we will define $B$-valued boxed convolution \frame{*}, as a binary
operation on $\Theta _{B}^{s}$ ; note that if $f,g\in \Theta _{B}^{s},$ then
we can always choose free $\{x_{1},...,x_{s}\}$ and $\{y_{1},...,y_{s}\}$ in
(some) NCPSpace over $B,$ $(A,\varphi ),$ such that

\begin{center}
$f=R_{x_{1},...,x_{s}}$ \ \ and \ \ $g=R_{y_{1},...,y_{s}}.$
\end{center}

\strut

\begin{definition}
(1) Let $s\in \Bbb{N}.$ Let $(f,g)\in \Theta _{B}^{s}\times \Theta _{B}^{s}.$
Define \frame{*}\thinspace \thinspace $:\Theta _{B}^{s}\times \Theta
_{B}^{s}\rightarrow \Theta _{B}^{s}$ by

\strut

\begin{center}
$\left( f,g\right) =\left( R_{x_{1},...,x_{s}},\,R_{y_{1},...,y_{s}}\right)
\longmapsto R_{x_{1},...,x_{s}}\,\,\frame{*}\,\,R_{y_{1},...,y_{s}}.$
\end{center}

\strut

Here, $\{x_{1},...,x_{s}\}$ and $\{y_{1},...,y_{s}\}$ are free in $%
(A,\varphi )$. Suppose that

\strut

\begin{center}
$coef_{i_{1},..,i_{n}}\left( R_{x_{1},...,x_{s}}\right)
=c^{(n)}(x_{i_{1}}\otimes b_{i_{2}}x_{i_{2}}\otimes ...\otimes
b_{i_{n}}x_{i_{n}})$
\end{center}

and

\begin{center}
$coef_{i_{1},...,i_{n}}(R_{y_{1},...,y_{s}})=c^{(n)}(y_{i_{1}}\otimes
b_{i_{2}}^{\prime }y_{i_{2}}\otimes ...\otimes b_{i_{n}}^{\prime
}y_{i_{n}}), $
\end{center}

\strut

for all $(i_{1},...,i_{n})\in \{1,...,s\}^{n},$ $n\in \Bbb{N},$ where $%
b_{i_{j}},b_{i_{n}}^{\prime }\in B$ arbitrary. Then

\strut

$coef_{i_{1},...,i_{n}}\left( R_{x_{1},...,x_{s}}\,\,\frame{*}%
\,\,R_{y_{1},...,y_{s}}\right) $

\strut

$=\underset{\pi \in NC(n)}{\sum }\left( \widehat{c_{x}}\oplus \widehat{c_{y}}%
\right) (\pi \cup Kr(\pi ))(x_{i_{1}}\otimes y_{i_{1}}\otimes
b_{i_{2}}x_{i_{2}}\otimes b_{i_{2}}^{\prime }y_{i_{2}}\otimes ...\otimes
b_{i_{n}}x_{i_{n}}\otimes b_{i_{n}}^{\prime }y_{i_{n}})$

$\strut $

$\overset{denote}{=}\underset{\pi \in NC(n)}{\sum }\left( k_{\pi }\oplus
k_{Kr(\pi )}\right) (x_{i_{1}},y_{i_{1}},...,x_{i_{n}}y_{i_{n}}),$

\strut \strut

where $\widehat{c_{x}}\oplus \widehat{c_{y}}=\widehat{c}\mid
_{A_{x}*_{B}A_{y}},$ $A_{x}=A\lg \left( \{x_{i}\}_{i=1}^{s},B\right) $ and $%
A_{y}=A\lg \left( \{y_{i}\}_{i=1}^{s},B\right) $ and where $\pi \cup Kr(\pi
) $ is an alternating union of partitions in $NC(2n)$
\end{definition}

\strut

\begin{proposition}
(See [12])\strut Let $(A,\varphi )$ be a NCPSpace over $B$ and let $%
x_{1},...,x_{s},y_{1},...,y_{s}\in (A,\varphi )$ be $B$-valued random
variables ($s\in \Bbb{N}$). If $\{x_{1},...,x_{s}\}$ and $%
\{y_{1},...,y_{s}\} $ are free in $(A,\varphi ),$ then we have

\strut

$k_{n}(x_{i_{1}}y_{i_{1}},...,x_{i_{n}}y_{i_{n}})$

\strut

$=\underset{\pi \in NC(n)}{\sum }\left( \widehat{c_{x}}\oplus \widehat{c_{y}}%
\right) (\pi \cup Kr(\pi ))(x_{i_{1}}\otimes y_{i_{1}}\otimes
b_{i_{2}}x_{i_{2}}\otimes y_{i_{2}}\otimes ...\otimes
b_{i_{n}}x_{i_{n}}\otimes y_{i_{n}})$

\strut

$\overset{denote}{=}\underset{\pi \in NC(n)}{\sum }\left( k_{\pi }\oplus
k_{Kr(\pi )}^{symm(1_{B})}\right)
(x_{i_{1}},y_{i_{1}},...,x_{i_{n}},y_{i_{n}}),$

\strut

for all $(i_{1},...,i_{n})\in \{1,...,s\}^{n},$ $\forall n\in \Bbb{N},$ $%
b_{i_{2}},...,b_{i_{n}}\in B,$ arbitrary, where $\widehat{c_{x}}\oplus 
\widehat{c_{y}}=\widehat{c}\mid _{A_{x}*_{B}A_{y}},$ $A_{x}=A\lg \left(
\{x_{i}\}_{i=1}^{s},B\right) $ and $A_{y}=A\lg \left(
\{y_{i}\}_{i=1}^{s},B\right) .$ \ $\square $
\end{proposition}

\strut

This shows that ;

\strut

\begin{corollary}
(See [12]) Under the same condition with the previous proposition,

\strut

\begin{center}
$R_{x_{1},...,x_{s}}\,\,\frame{*}\,%
\,R_{y_{1},...,y_{s}}^{t}=R_{x_{1}y_{1},...,x_{s}y_{s}}.$
\end{center}

$\square $
\end{corollary}

\strut

Notice that, in general, unless $b_{i_{2}}^{\prime }=...=b_{i_{n}}^{\prime
}=1_{B}$ in $B,$

\strut

\begin{center}
$R_{x_{1},...,x_{s}}\,\,\frame{*}\,\,R_{y_{1},...,y_{s}}\neq
R_{x_{1}y_{1},...,x_{s}y_{s}}.$
\end{center}

\strut

However, as we can see above,

\strut

\begin{center}
$R_{x_{1},...,x_{s}}\,\,\frame{*}\,%
\,R_{y_{1},...,y_{s}}^{t}=R_{x_{1}y_{1},...,x_{s}y_{s}}$
\end{center}

and

\begin{center}
$R_{x_{1},...,x_{s}}^{t}\,\,\frame{*}\,%
\,R_{y_{1},...,y_{s}}^{t}=R_{x_{1}y_{1},...,x_{s}y_{s}}^{t},$
\end{center}

\strut

where $\{x_{1},...,x_{s}\}$ and $\{y_{1},...,y_{s}\}$ are free over $B.$
Over $B=\Bbb{C},$ the last equation is proved by Nica and Speicher in [2]
and [3]. Actually, their R-transforms (over $\Bbb{C}$) is our trivial
R-transforms (over $\Bbb{C}$). Is it possible to find a trivial R-transform $%
R_{x_{1}^{\prime },...,x_{s}^{\prime }}^{t},$ for the given R-transform $%
R_{x_{1},...,x_{s}},$ where $x_{1},...,x_{s}\in (A,\varphi )$ are $B$-valued
random variables and $x_{1}^{\prime },...,x_{s}^{\prime }$ are $B$-valued
random variables (which are not necessarily contained in $(A,\varphi )$%
)\thinspace ? If possible, to study amalgamated R-transforms, we can only
deal with trivial R-transforms. Let's denote a set of all tritial
R-transforms by $S^{(1_{B})}\,\mathcal{R}_{B}^{s}.$ If we have a positive
answer of above question, then $\mathcal{R}_{B}^{s}=S^{(1_{B})}\,\mathcal{R}%
_{B}^{s}$. In that case, we can deal with $B$-valued boxed convolution more
freely. In the following chapter, we will consider this possibility.

\strut

\strut

\strut

\strut

\subsection{Scalar-Valued Case (in the sense of Nica and Speicher)}

\strut

\strut

\strut

Suppose that $B=\Bbb{C}.$ Then a NCPSpace over $B,$ $(A,\varphi )$, All
information of $B$-valued (scalar-valued) distributions are contained in
corresponding\ trivial R-transforms (See [32] and [33]). So, WLOG, we can
define an R-transform of random variable(s) as a trivial R-transform of the
random variable(s). This is just an (scalar-valued) R-transform in the sense
of Nica (See [2]). In [32] and [33], we showed that if $B$ is an arbitrary
algebra such that $B=C_{A}(B),$ then amalgamated freeness can be
characterized by

\strut

Amalgamated Freeness $\overset{def}{\Leftrightarrow }$ All mixed cumulants
vanish

\begin{center}
$\Leftrightarrow $ All mixed\textbf{\ trivial }cumulants vanish.
\end{center}

\strut

Since $C_{A}(\Bbb{C})=\Bbb{C},$ for any algebra $A,$ if $B=\Bbb{C},$ then
the freeness is characterized by trivial cumulants.

\strut

\begin{definition}
Let $(A,\varphi )$ be a noncommutative probability space over $\Bbb{C}$
(shortly NCPSpace) with a linear functional $\varphi :A\rightarrow \Bbb{C}$
and let $x_{1},...,x_{s}\in (A,\varphi )$ be (scalar-valued) random
variables ($s\in \Bbb{N}$). Define a moment series and R-transform of $%
x_{1},...,x_{s}$ by $M_{x_{1},...,x_{s}}^{t}$ and $R_{x_{1},...,x_{s}}^{t},$
respectively, where $M_{x_{1},...,x_{s}}^{t}$ and $R_{x_{1},...,x_{s}}^{t}$
are defined in Section 1.2. By $M_{x_{1},...,x_{s}}$ and $%
R_{x_{1},...,x_{s}},$ we will denote the moment series and the R-transform,
respectively. Also, if we mention cumulants, then they are trivial
cumulants, just denoted by $k_{n}(\cdot )$, for $n\in \Bbb{N}.$
\end{definition}

\strut

\begin{remark}
By definition, all facts, which we considered in the previous sections,
holds true, for new $M_{x_{1},...,x_{s}}$ and $R_{x_{1},...,x_{s}}.$
Moreover, if $\{x_{1},...,x_{s}\}$ and $\{y_{1},...,y_{s}\}$ are free
subsets of random variables, then

\strut

\begin{center}
$R_{x_{1}y_{1},...,x_{s}y_{s}}=R_{x_{1},...,x_{s}}\,\,\frame{*}_{\Bbb{C}%
}\,\,R_{y_{1},...,y_{s}}$,
\end{center}

\strut

where \frame{*}$_{\Bbb{C}}$ is the $B=\Bbb{C}$-valued boxed convolution.
This is the boxed convolution \frame{*} in the sense of Nica (See [2]).
Notice that all random variables are central over $\Bbb{C}.$ So, we have that

\strut

$coef_{i_{1},...,i_{n}}\left( R_{x_{1},...,x_{s}}\,\,\frame{*}%
\,\,R_{y_{1},...,y_{s}}\right) $

\begin{center}
$=\underset{\pi \in NC(n)}{\sum }\left( k_{\pi
}(x_{i_{1}},...,x_{i_{n}})\right) \left( k_{Kr(\pi
)}(y_{i_{1}},...,y_{i_{n}})\right) .$
\end{center}

\strut

Again, we can see that this is just a definition of boxed convolution in the
sense of Nica and Speicher (See [2] and [3]).
\end{remark}

\strut

Form now, by $(A,\varphi ),$ we will denote a NCPSpace.

\strut

\strut

\strut

\strut

\strut

\strut

\section{Toeplitz Matricial Probability Theory}

\strut

\strut

\strut

\strut

\subsection{Toeplitz Matricial Probability Spaces over Toeplitz Matricial
Algebras}

\strut

\strut

\strut

In this section, we will consider the Toeplitz Matricial probability space
over a Toeplitz Matricial Algebra $\mathcal{C}^{n}\cong \Bbb{C}^{n},$ for $%
n\in \Bbb{N}.$ We define a Toeplitz Matricial algebra $\mathcal{C}^{n}$ by
an algebra $\Bbb{C}^{n}$ with following addition and multiplication ;

\strut

\begin{center}
$(\alpha _{1},...,\alpha _{n})+(\alpha _{1}^{\prime },...,\alpha
_{n}^{\prime })=(\alpha _{1}+\alpha _{1}^{\prime },...,\alpha _{n}+\alpha
_{n}^{\prime })$
\end{center}

and

\begin{center}
$(\alpha _{1},...,\alpha _{n})\cdot (\alpha _{1}^{\prime },...,\alpha
_{n}^{\prime })=(\alpha _{1}\alpha _{1}^{\prime },\,\sum_{k=1}^{2}\alpha
_{k}\alpha _{(2+1)-k}^{\prime },\,...,\,\sum_{k=1}^{n}\alpha _{k}\alpha
_{(n+1)-k}^{\prime }),$
\end{center}

\strut

respectively. Then $\mathcal{C}^{n}=\left( \Bbb{C}^{n},+,\,\cdot \right) $
is an algebra. And this algebra is called an $n$-th Toeplitz Matricial
algebra. i.e isometrically, this algebra $\mathcal{C}^{n}$ can be understood
as an algebra generated by all Toeplitz matrices in $M_{n}(\Bbb{C})$ which
have the form of

\strut

\begin{center}
$(\alpha _{1},...,\alpha _{n})=\left( 
\begin{array}{lllll}
\alpha _{1} & \alpha _{2} & \alpha _{3} & \cdots & \alpha _{n} \\ 
& \alpha _{1} & \alpha _{2} & \ddots & \vdots \\ 
&  & \alpha _{1} & \ddots & \alpha _{3} \\ 
&  &  & \ddots & \alpha _{2} \\ 
O &  &  &  & \alpha _{1}
\end{array}
\right) \in M_{n}(\Bbb{C}),$
\end{center}

\strut

for all $(\alpha _{1},...,\alpha _{n})\in \mathcal{C}^{n},$ $n\in \Bbb{N}.$

\strut

For exmaple, if $n=4,$ then we have that

\strut

$(\alpha _{1},\alpha _{2},\alpha _{3},\alpha _{4})+(\alpha _{1}^{\prime
},\alpha _{2}^{\prime },\alpha _{3}^{\prime },\alpha _{4}^{\prime })$

$\ \ \cong \left( 
\begin{array}{llll}
\alpha _{1} & \alpha _{2} & \alpha _{3} & \alpha _{4} \\ 
0 & \alpha _{1} & \alpha _{2} & \alpha _{3} \\ 
0 & 0 & \alpha _{1} & \alpha _{2} \\ 
0 & 0 & 0 & \alpha _{1}
\end{array}
\right) +\left( 
\begin{array}{llll}
\alpha _{1}^{\prime } & \alpha _{2}^{\prime } & \alpha _{3}^{\prime } & 
\alpha _{4}^{\prime } \\ 
0 & \alpha _{1}^{\prime } & \alpha _{2}^{\prime } & \alpha _{3}^{\prime } \\ 
0 & 0 & \alpha _{1}^{\prime } & \alpha _{2}^{\prime } \\ 
0 & 0 & 0 & \alpha _{1}^{\prime }
\end{array}
\right) $

$\ \ =\left( 
\begin{array}{cccc}
\alpha _{1}+\alpha _{1}^{\prime } & \alpha _{2}+\alpha _{2}^{\prime } & 
\alpha _{3}+\alpha _{3}^{\prime } & \alpha _{4}+\alpha _{4}^{\prime } \\ 
0 & \alpha _{1}+\alpha _{1}^{\prime } & \alpha _{2}+\alpha _{2}^{\prime } & 
\alpha _{3}+\alpha _{3}^{\prime } \\ 
0 & 0 & \alpha _{1}+\alpha _{1}^{\prime } & \alpha _{2}+\alpha _{2}^{\prime }
\\ 
0 & 0 & 0 & \alpha _{1}+\alpha _{1}^{\prime }
\end{array}
\right) $

\strut

$\ \ \cong (\alpha _{1}+\alpha _{1}^{\prime },\,\alpha _{2}+\alpha
_{2}^{\prime },\,\alpha _{3}+\alpha _{3}^{\prime },\,\alpha _{4}+\alpha
_{4}^{\prime })$

\strut

and

\strut

$(\alpha _{1},\alpha _{2},\alpha _{3},\alpha _{4})\cdot (\alpha _{1}^{\prime
},\alpha _{2}^{\prime },\alpha _{3}^{\prime },\alpha _{4}^{\prime })$

\strut

$\cong \left( 
\begin{array}{llll}
\alpha _{1} & \alpha _{2} & \alpha _{3} & \alpha _{4} \\ 
0 & \alpha _{1} & \alpha _{2} & \alpha _{3} \\ 
0 & 0 & \alpha _{1} & \alpha _{2} \\ 
0 & 0 & 0 & \alpha _{1}
\end{array}
\right) \left( 
\begin{array}{llll}
\alpha _{1}^{\prime } & \alpha _{2}^{\prime } & \alpha _{3}^{\prime } & 
\alpha _{4}^{\prime } \\ 
0 & \alpha _{1}^{\prime } & \alpha _{2}^{\prime } & \alpha _{3}^{\prime } \\ 
0 & 0 & \alpha _{1}^{\prime } & \alpha _{2}^{\prime } \\ 
0 & 0 & 0 & \alpha _{1}^{\prime }
\end{array}
\right) $

$=\left( 
\begin{array}{cccc}
\alpha _{1}\alpha _{1}^{\prime } & \alpha _{1}\alpha _{2}^{\prime }+\alpha
_{2}\alpha _{1}^{\prime } & \alpha _{1}\alpha _{3}^{\prime }+\alpha
_{2}\alpha _{2}^{\prime }+\alpha _{3}\alpha _{1}^{\prime } & \alpha
_{1}\alpha _{4}^{\prime }+\alpha _{2}\alpha _{3}^{\prime }+\alpha _{3}\alpha
_{2}^{\prime }+\alpha _{4}\alpha _{1}^{\prime } \\ 
0 & \alpha _{1}\alpha _{1}^{\prime } & \alpha _{1}\alpha _{2}^{\prime
}+\alpha _{2}\alpha _{1}^{\prime } & \alpha _{1}\alpha _{3}^{\prime }+\alpha
_{2}\alpha _{2}^{\prime }+\alpha _{3}\alpha _{1}^{\prime } \\ 
0 & 0 & \alpha _{1}\alpha _{1}^{\prime } & \alpha _{1}\alpha _{2}^{\prime
}+\alpha _{2}\alpha _{1}^{\prime } \\ 
0 & 0 & 0 & \alpha _{1}\alpha _{1}^{\prime }
\end{array}
\right) $

\strut

$\cong (\alpha _{1}\alpha _{1}^{\prime },\,\alpha _{1}\alpha _{2}^{\prime
}+\alpha _{2}\alpha _{1}^{\prime },\,\alpha _{1}\alpha _{3}^{\prime }+\alpha
_{2}\alpha _{2}^{\prime }+\alpha _{3}\alpha _{1}^{\prime },\,\alpha
_{1}\alpha _{4}^{\prime }+\alpha _{2}\alpha _{3}^{\prime }+\alpha _{3}\alpha
_{2}^{\prime }+\alpha _{4}\alpha _{1}^{\prime }),$

\strut

for all $(\alpha _{1},...,\alpha _{4}),\,(\alpha _{1}^{\prime },...,\alpha
_{4}^{\prime })\in \mathcal{C}^{4}.$ Notice that $(1,0,...,0)\in \mathcal{C}%
^{n}$ is the multiplication-identity on $\mathcal{C}^{n}.$

\strut

We will observe the R-transform calculus in a noncommutative probability
space over a Toeplitz matricial algebra $\mathcal{C}^{n}.$

\strut

\begin{quote}
\frame{\textbf{Notation}} From now, for any $n\in \Bbb{N},$ by little abuse
of notation, we will denote $\mathcal{C}^{n}$ by $B,$ if there is no
confusion. $\square $
\end{quote}

\strut

Let $(A,\varphi )$ be a NCPSpace in the sense of Section 1.3, where $A$ is
an algebra and $\varphi :A\rightarrow \Bbb{C}$ is a linear functional. Then
we can define a Toeplitz matricial algebra $T_{n}$ over a Toeplitz matricial
algebra $B=\mathcal{C}^{n},$ induced by $(A,\varphi )$ ;

\strut

\begin{definition}
Fix $n\in \Bbb{N}.$ Let $B=\mathcal{C}^{n}$ be a Toeplitz matricial algebra
and let $(A,\varphi )$ be a NCPSpace in the sense of Section 1.3, with a
linear functional $\varphi :A\rightarrow \Bbb{C}.$ Define a Toeplitz
Matricial Probability Space over $B,$ $(T_{n},E_{n}),$ induced by $%
(A,\varphi ),$ by

\strut

(i) $T_{n}$ is an algebra over $B$ (i.e $1_{B}=(1,0,...,0)\in B\subset T_{n}$%
) such that

\strut

\begin{center}
$T_{n}=A\lg \left( \{(a_{1},...,a_{n}):\forall \,a_{j}\in (A,\varphi
)\}\right) $
\end{center}

with

\strut

\begin{center}
$(a_{1},...,a_{n})+(a_{1}^{\prime },...,a_{n}^{\prime
})=(a_{1}+a_{1}^{\prime },...,a_{n}+a_{n}^{\prime })$
\end{center}

and

\strut

\begin{center}
$(a_{1},...,a_{n})\cdot (a_{1}^{\prime },...,a_{n}^{\prime
})=(a_{1}a_{1}^{\prime },\,\sum_{k=1}^{2}a_{k}a_{(2+1)-k}^{\prime
},...,\sum_{k=1}^{n}a_{k}a_{(n+1)-k}^{\prime }),$
\end{center}

\strut

for all $(a_{1},...,a_{n}),$ $(a_{1}^{\prime },...,a_{n}^{\prime })\in
T_{n}. $

\strut

(ii) $E:T_{n}\rightarrow B$ is a $B$-functional

\strut

\begin{center}
$E\left( (a_{1},...,a_{n})\right) =\left( \varphi (a_{1}),...,\varphi
(a_{n})\right) .$
\end{center}

\strut

Without confusion, we will denote $(T_{n},E_{n})$ by $(T,E).$
\end{definition}

\strut

\begin{lemma}
Let $n\in \Bbb{N}.$ A Toeplitz matricial probability space, $\left(
T,E\right) ,$ over $B=\mathcal{C}^{n},$ is, indeed, a NCPSPace over $B.$
\end{lemma}

\strut

\begin{proof}
Since $T$ is an algebra over $B,$ by definition, it suffices to show that a
map $E:T\rightarrow B$ is a $B$-functional.

\strut

(i) \ Let $(\alpha _{1},...,\alpha _{n})\in \mathcal{C}^{n}=B$. Then

\strut

\begin{center}
$E\left( (\alpha _{1},...,\alpha _{n})\right) =\left( \varphi (\alpha
_{1}),...,\varphi (\alpha _{n})\right) =(\alpha _{1},...,\alpha _{n}),$
\end{center}

\strut

since $\varphi $ is a linear functional.

\strut

(ii) Let $(\alpha _{1},...,\alpha _{n}),(\alpha _{1}^{\prime },...,\alpha
_{n}^{\prime })\in B$ and $(a_{1},...,a_{n})\in T.$ Then

\strut

$E\left( (\alpha _{1},...,\alpha _{n})(a_{1},...,a_{n})(\alpha _{1}^{\prime
},...,\alpha _{n}^{\prime })\right) $

\strut

$\equiv E\left( \left( 
\begin{array}{llll}
\alpha _{1} & \alpha _{2} & \cdots & \alpha _{n} \\ 
& \alpha _{1} & \ddots & \vdots \\ 
&  & \ddots & \alpha _{2} \\ 
O &  &  & \alpha _{1}
\end{array}
\right) \left( 
\begin{array}{llll}
a_{1} & a_{2} & \cdots & a_{n} \\ 
& a_{1} & \ddots & \vdots \\ 
&  & \ddots & a_{2} \\ 
O &  &  & a_{1}
\end{array}
\right) \left( 
\begin{array}{llll}
\alpha _{1}^{\prime } & \alpha _{2}^{\prime } & \cdots & \alpha _{n}^{\prime
} \\ 
& \alpha _{1}^{\prime } & \ddots & \vdots \\ 
&  & \ddots & \alpha _{2}^{\prime } \\ 
&  &  & \alpha _{1}^{\prime }
\end{array}
\right) \right) $

\strut

$=\left( 
\begin{array}{llll}
\alpha _{1} & \alpha _{2} & \cdots & \alpha _{n} \\ 
& \alpha _{1} & \ddots & \vdots \\ 
&  & \ddots & \alpha _{2} \\ 
O &  &  & \alpha _{1}
\end{array}
\right) \cdot \left( 
\begin{array}{llll}
\varphi (a_{1}) & \varphi (a_{2}) & \cdots & \varphi (a_{n}) \\ 
& \varphi (a_{1}) & \ddots & \vdots \\ 
&  & \ddots & \varphi (a_{2}) \\ 
&  &  & \varphi (a_{1})
\end{array}
\right) \cdot \left( 
\begin{array}{llll}
\alpha _{1}^{\prime } & \alpha _{2}^{\prime } & \cdots & \alpha _{n}^{\prime
} \\ 
& \alpha _{1}^{\prime } & \ddots & \vdots \\ 
&  & \ddots & \alpha _{2}^{\prime } \\ 
&  &  & \alpha _{1}^{\prime }
\end{array}
\right) $

\strut

$\equiv (\alpha _{1},...,\alpha _{n})\cdot E\left( (a_{1},...,a_{n})\right)
\cdot (\alpha _{1}^{\prime },...,\alpha _{n}^{\prime }).$

\strut

By (i) and (ii), we can conclude that $E$ is a $B=\mathcal{C}^{n}$%
-functional.
\end{proof}

\strut

By the previous lemma, we have a new NCPSpace over $B=\mathcal{C}^{n},$ $%
(T,E),$ in the sense of Section 1.1. On this framework, we can do the
amalgamated R-transform Calculus, like in Section 1.2.

\strut

\begin{lemma}
Fix $n\in \Bbb{N}.$ Let $B=\mathcal{C}^{n}$ be a Toeplitz matricial algebra
and let $(A,\varphi )$ be a NCPSpace. Let $(T,E)$ be a Toeplitz matricial
probability space over $B,$ induced by $(A,\varphi ).$ Then $C_{T}(B)=B,$
where $C_{T}(B)$ is the centralizer of $B$ in $T,$ i.e

\strut

\begin{center}
$C_{T}(B)=\{b\in B:bt=ta,$ for all $t\in T\}.$
\end{center}
\end{lemma}

\strut

\begin{proof}
Let $(\alpha _{1},...,\alpha _{n})\in \mathcal{C}^{n}=B.$ Then $(\alpha
_{1},...,\alpha _{n})\in C_{T}(B).$ i.e

\strut

$(\alpha _{1},...,\alpha _{n})\cdot (a_{1},...,a_{n})$

$\ \ \ =\left( \alpha _{1}a_{1},\sum_{k=1}^{2}\alpha
_{k}a_{(2+1)-k},...,\sum_{k=1}^{n}\alpha _{k}a_{(n+1)-k}\right) $

$\ \ \ =\left( a_{1}\alpha _{1},\sum_{k=1}^{2}a_{(2+1)-k}\alpha
_{k},...,\sum_{k=1}^{n}a_{(n+1)-k}\alpha _{k}\right) $

$\ \ \ =\left( a_{1}\alpha _{1},\sum_{j=1}^{2}a_{j}\alpha
_{(2+1)-j},...,\sum_{j=1}^{n}a_{j}\alpha _{(n+1)-j}\right) $

$\ \ \ =(a_{1},...,a_{n})\cdot (\alpha _{1},...,\alpha _{n}),$

\strut

for all $(a_{1},...,a_{n})\in T.$

\strut

Trivially, $C_{T}(B)\subseteq B,$ by the very definition of $C_{T}(B).$
Therefore,

\strut

\begin{center}
$C_{T}(B)=B.$
\end{center}
\end{proof}

\strut

Fix $N\in \Bbb{N}.$ Like in Section 1.3, since $C_{T}(B)=B,$ trivial
R-transforms of $B=\mathcal{C}^{N}$-valued random variables contain all
information of operator-valued distributions of those $B$-valued random
variables (Also see [32] and [33]). Hence, it is sufficient to consider
trivial R-transforms of $B$-valued random variables to study operator valued
random variables, in $\sum_{B}^{s},$ where $B=\mathcal{C}^{N}$ and $s\in 
\Bbb{N}.$ Remark that to study such trivial R-transforms, we need to define
suitable cumulants $(K_{n})_{n=1}^{\infty },$ as a cumulant multiplicative
bimodule map induced by a operator-valued (moment multiplicative)
conditional expectation $E.$ Recall that

\strut

\begin{center}
$E\left( (a_{1},...,a_{n})\right) \overset{def}{=}\left( \varphi
(a_{1}),...,\varphi (a_{n})\right) \in \mathcal{C}^{n},$
\end{center}

\strut

for all $(a_{1},...,a_{N})\in T.$ Denote cumulants induced by $\varphi ,$
with respect to $(A,\varphi ),$ by $(k_{n})_{n=1}^{\infty }.$ i,e if $%
x_{1},...,x_{s}\in (A,\varphi )$ are random variables, then

\strut

\begin{center}
$k_{n}\left( x_{i_{1}},.....,x_{i_{n}}\right) \overset{def}{=}\underset{\pi
\in NC(n)}{\sum }\widehat{\varphi }(\pi )\left( x_{i_{1}}\otimes
.....\otimes x_{i_{n}}\right) \mu (\pi ,1_{n})\in \Bbb{C},$
\end{center}

\strut

for all $(i_{1},...,i_{n})\in \{1,...,s\}^{n},$ $n\in \Bbb{N}$, by the
M\"{o}bius inversion. \strut

\strut

Suppose that $(a_{1}^{(1)},...,a_{N}^{(1)})$,...., $%
(a_{1}^{(s)},...,a_{N}^{(s)})\in T_{N},$ for the fixed $N\in \Bbb{N}.$
Observe the following relation, for $(i_{1},...,i_{n})\in \{1,...,s\}^{n},$ $%
n\in \Bbb{N}$ ;

\strut

$%
(a_{1}^{(i_{1})},...,a_{N}^{(i_{1})})(a_{1}^{(i_{2})},...,a_{N}^{(i_{2})})=(a_{1}^{(i_{1})}a_{1}^{(i_{2})},\sum_{k=1}^{2}a_{k}^{(i_{1})}a_{(2+1)-k}^{(i_{2})},...,\sum_{k=1}^{N}a_{k}^{(i_{1})}a_{(N+1)-k}^{(i_{2})}). 
$

\strut

Put\strut

\begin{center}
$P_{j}^{(i_{1},i_{2})}=\sum_{k=1}^{j}a_{k}^{(i_{1})}a_{(j+1)-k}^{(i_{2})}\in
(A,\varphi )$ \ \ \ ( $j=1,...,N$).
\end{center}

i.e

\begin{center}
$(a_{1}^{(i_{1})}a_{1}^{(i_{2})},%
\sum_{k=1}^{2}a_{k}^{(i_{1})}a_{(2+1)-k}^{(i_{2})},...,%
\sum_{k=1}^{N}a_{k}^{(i_{1})}a_{(N+1)-k}^{(i_{2})})=\left(
P_{1}^{(i_{1},i_{2})},...,P_{N}^{(i_{1},i_{2})}\right) .$
\end{center}

\strut

Then

$%
(a_{1}^{(i_{1})},...,a_{N}^{(i_{1})})(a_{1}^{(i_{2})},...,a_{N}^{(i_{2})})(a_{1}^{(i_{3})},...,a_{N}^{(i_{3})}) 
$

$\ \ =\left( P_{1}^{(i_{1},i_{2})},...,P_{N}^{(i_{1},i_{2})}\right)
(a_{1}^{(i_{3})},...,a_{N}^{(i_{3})})$

$\ \ =\left(
P_{1}^{(i_{1},i_{2})}a_{1}^{(i_{3})},%
\sum_{k=1}^{2}P_{k}^{(i_{1},i_{2})}a_{(2+1)-k}^{(i_{3})},...,%
\sum_{k=1}^{N}P_{N}^{(i_{1},i_{2})}a_{(N+1)-k}^{(i_{3})}\right) .$

\strut

Similarly, put

\strut

\begin{center}
$P_{j}^{(i_{1},i_{2},i_{3})}=%
\sum_{k=1}^{j}P_{k}^{(i_{1},i_{2})}a_{(j+1)-k}^{(i_{3})}$ \ \ \ ($j=1,...,N$%
).
\end{center}

i.e

$\left(
P_{1}^{(i_{1},i_{2})}a_{1}^{(i_{3})},%
\sum_{k=1}^{2}P_{k}^{(i_{1},i_{2})}a_{(2+1)-k}^{(i_{3})},...,%
\sum_{k=1}^{N}P_{N}^{(i_{1},i_{2})}a_{(N+1)-k}^{(i_{3})}\right) $

\begin{center}
$=\left( P_{1}^{(i_{1},i_{2},i_{3})},...,P_{N}^{(i_{1},i_{2},i_{3})}\right)
. $
\end{center}

\strut

Inductively, we can define

\strut

\begin{center}
$P_{j}^{(i_{1},...,i_{n-1},i_{n})}=%
\sum_{k=1}^{j}P_{k}^{(i_{1},...,i_{n-1})}a_{(j+1)-k}^{(i_{n})}$ \ \ ($%
j=1,...,N$).
\end{center}

i.e

\strut

\begin{center}
$(a_{1}^{(i_{1})},...,a_{N}^{(i_{1})})\cdot \cdot \cdot
(a_{1}^{(i_{n})},...,a_{N}^{(i_{n})})=\left(
P_{1}^{(i_{1},..,i_{n})},...,P_{N}^{(i_{1},...,i_{n})}\right) .$
\end{center}

\strut

\begin{definition}
Let $B=\mathcal{C}^{N}$ be a Toeplitz matricial algebra and let $(A,\varphi
) $ be a NCPSpace with its linear functional $\varphi :A\rightarrow \Bbb{C}.$
Let $(T_{N},E_{N})\equiv (T,E)$ be a Toeplitz matricial probability space
over $B$ with its $B$-functional $E:T\rightarrow B$ and let $%
(a_{1}^{(1)},...,a_{N}^{(1)}),...,(a_{1}^{(s)},...,a_{N}^{(s)})\in (T,E)$ be 
$B$-valued random variables ($s\in \Bbb{N}$).

\strut

(1) For $(i_{1},...,i_{n})\in \{1,...,s\}^{n},$ $n\in \Bbb{N},$ define

\strut \strut

\begin{center}
$P_{j}^{(i_{1})}=a_{j}^{(i_{1})}\in (A,\varphi ),$ \ $j=1,...,N,$
\end{center}

\strut

\begin{center}
$P_{j}^{(i_{1},i_{2})}=\sum_{k=1}^{j}a_{k}^{(i_{1})}a_{(j+1)-k}^{(i_{2})},$
\ $j=1,...,N$
\end{center}

and

\strut

\begin{center}
$P_{j}^{(i_{1},...,i_{n-1},i_{n})}=%
\sum_{k=1}^{j}P_{k}^{(i_{1},...,i_{n-1})}a_{(j+1)-k}^{(i_{n})},$ \ $%
j=1,...,N.$
\end{center}

\strut

Then

\begin{center}
$(a_{1}^{(i_{1})},...,a_{N}^{(i_{1})})\cdot \cdot \cdot
(a_{1}^{(i_{n})},...,a_{N}^{(i_{n})})=\left(
P_{1}^{(i_{1},...,i_{n})},...,P_{N}^{(i_{1},...,i_{n})}\right) .$
\end{center}

\strut

(2) Define a $(i_{1},...,i_{n})$-th moment of $%
(a_{1}^{(1)},...,a_{N}^{(1)}),...,(a_{1}^{(s)},...,a_{N}^{(s)})$ by

\strut

\begin{center}
$E\left( (a_{1}^{(i_{1})},...,a_{N}^{(i_{1})})\cdot \cdot \cdot
(a_{1}^{(i_{n})},...,a_{N}^{(i_{n})})\right) =\left( \varphi
(P_{1}^{(i_{1},...,i_{n})}),...,\varphi (P_{N}^{(i_{1},...,i_{n})})\right) ,$
\end{center}

\strut

for all $(i_{1},...,i_{n})\in \{1,...,s\}^{n},$ $n\in \Bbb{N}.$
\end{definition}

\strut

Notice that

\strut

\begin{center}
$\varphi \left( P_{1}^{(i_{1},...,i_{n})}\right) =\varphi \left(
a_{1}^{(i_{1})}...a_{1}^{(i_{n})}\right) $
\end{center}

\strut

is the $(i_{1},...,i_{n})$-th moment of $a_{1}^{(1)},...,a_{1}^{(s)}.$ As we
have seen before (and in [32] and [33]), all \textbf{trivial} moments
contain all information about moments, since $C_{T}(B)=B.$ Hence WLOG, we
can only consider the trivial moments of $(a_{1}^{(1)},...,a_{N}^{(1)}),...,%
\,(a_{1}^{(s)},...,a_{N}^{(s)}).$

\strut

\begin{example}
Suppose that $N=2$ and fix \ $B=\mathcal{C}^{2}.$ Let $(T_{2},E_{2})$ be a
Toeplitz matricial probability space over $B,$ induced by a NCPSpace $%
(A,\varphi ).$ Let $(a_{1}^{(1)},a_{2}^{(1)}),$ $(a_{1}^{(2)},a_{2}^{(2)}),$ 
$(a_{1}^{(3)},a_{2}^{(3)})\in (T_{2},E_{2})$ be $B$-valued random variables.
Then

\strut

$E\left( (a_{1}^{(1)},a_{2}^{(1)})\cdot (a_{1}^{(2)},a_{2}^{(2)})\cdot
(a_{1}^{(3)},a_{2}^{(3)})\right) $

\strut

$\ \ \equiv E\left( \left( 
\begin{array}{ll}
a_{1}^{(1)} & a_{2}^{(1)} \\ 
O & a_{1}^{(1)}
\end{array}
\right) \left( 
\begin{array}{ll}
a_{1}^{(2)} & a_{2}^{(2)} \\ 
O & a_{1}^{(2)}
\end{array}
\right) \left( 
\begin{array}{ll}
a_{1}^{(3)} & a_{2}^{(3)} \\ 
O & a_{1}^{(3)}
\end{array}
\right) \right) $

\strut

where $O=0_{B}=0_{A}\in A.$

\strut

$\ \ =E\left( \left( 
\begin{array}{lll}
a_{1}^{(1)}a_{1}^{(2)}a_{1}^{(3)} &  & 
a_{1}^{(1)}a_{1}^{(2)}a_{2}^{(3)}+a_{1}^{(1)}a_{2}^{(2)}a_{1}^{(3)}+a_{2}^{(1)}a_{1}^{(2)}a_{1}^{(3)}
\\ 
&  &  \\ 
O &  & \text{ \ \ \ \ \ \ \ \ \ \ \ \ \ \ \ \ \ \ \ \ \ }%
a_{1}^{(1)}a_{1}^{(2)}a_{1}^{(3)}
\end{array}
\right) \right) $\strut

\strut

$\ \ \equiv E\left(
(a_{1}^{(1)}a_{1}^{(2)}a_{1}^{(3)},\,\,\,%
\,a_{1}^{(1)}a_{1}^{(2)}a_{2}^{(3)}+a_{1}^{(1)}a_{2}^{(2)}a_{1}^{(3)}+a_{2}^{(1)}a_{1}^{(2)}a_{1}^{(3)})\right) 
$

\strut

$\ \ =\left( \varphi (a_{1}^{(1)}a_{1}^{(2)}a_{1}^{(3)}),\,\,\,\varphi
(a_{1}^{(1)}a_{1}^{(2)}a_{2}^{(3)}+a_{1}^{(1)}a_{2}^{(2)}a_{1}^{(3)}+a_{2}^{(1)}a_{1}^{(2)}a_{1}^{(3)})\right) . 
$

\strut

We have that

\strut

\begin{center}
$P_{1}^{(1,2,3)}=a_{1}^{(1)}a_{1}^{(2)}a_{1}^{(3)}$
\end{center}

and

\begin{center}
$%
P_{2}^{(1,2,3)}=a_{1}^{(1)}a_{1}^{(2)}a_{2}^{(3)}+a_{1}^{(1)}a_{2}^{(2)}a_{1}^{(3)}+a_{2}^{(1)}a_{1}^{(2)}a_{1}^{(3)}. 
$
\end{center}
\end{example}

\strut

\strut

\strut

\strut

\strut

\strut

\subsection{Toeplitz Matricial Cumulants}

\strut

\strut

\strut

\strut

Throughout this section, we will fix $N\in \Bbb{N}.$ Let $B=\mathcal{C}^{N}$
be a Toeplitz matricial algebra and let $(A,\varphi )$ be a NCPSpace, with
its linear functional $\varphi :A\rightarrow \Bbb{C}.$ Let $%
(T_{N},E_{N})\equiv (T,E)$ be a Toeplitz matricial probability space over $B$
with its $B$-functional $E:T\rightarrow B,$

\strut

\begin{center}
$E\left( (a_{1},...,a_{N})\right) =\left( \varphi (a_{1}),...,\varphi
(a_{N})\right) \in B,$
\end{center}

\strut

for all $(a_{1},...,a_{N})\in T.$\strut Let $%
(a_{1}^{(1)},...,a_{N}^{(1)}),...,(a_{1}^{(s)},...,a_{N}^{(s)})\in (T,E)$ be 
$B$-valued random variables ($s\in \Bbb{N}$). Then we can conpute the $%
(i_{1},...,i_{n})$-th moment of them by

\strut

$E\left(
(a_{1}^{(i_{1})},...,a_{N}^{(i_{1})})...(a_{1}^{(i_{n})},...,a_{N}^{(i_{n})})\right) 
$

$\ \ \ \ \ \ \ \ \ \ \ \ \ \ \ \ \ =E\left(
(P_{1}^{(i_{1},...,i_{n})},...,P_{N}^{(i_{1},...,i_{n})})\right) $

$\ \ \ \ \ \ \ \ \ \ \ \ \ \ \ \ \ =\left( \varphi
(P_{1}^{(i_{1},...,i_{n})}),....,\varphi (P_{N}^{(i_{1},...,i_{n})})\right)
, $

\strut

where

\strut

\begin{center}
$P_{j}^{(i_{1})}=a_{j}^{(i_{1})}\in (A,\varphi ),$
\end{center}

\strut

\begin{center}
$P_{j}^{(i_{1},i_{2})}=\sum_{k=1}^{j}a_{k}^{(i_{1})}a_{(j+1)-k}^{(i_{2})}\in
(A,\varphi ),$
\end{center}

and

\begin{center}
$P_{j}^{(i_{1},...,i_{n-1},i_{n})}=%
\sum_{k=1}^{j}P_{k}^{(i_{1},...,i_{n-1})}a_{(j+1)-k}^{(i_{n})}\in (A,\varphi
),$
\end{center}

\strut

for all $j=1,...,N.$

\strut

So, we can define a moment series of $%
A_{1}=(a_{1}^{(1)},...,a_{N}^{(1)}),...,A_{s}=(a_{1}^{(s)},...,a_{N}^{(s)})$
denoted by $M_{A_{1},...,A_{s}},$ by

\strut

$%
M_{A_{1},...,A_{s}}(z_{1},...,z_{s})=M_{A_{1},...,A_{s}}^{t}(z_{1},...,z_{s}) 
$

$\ \ \ \ \ \ =\sum_{n=1}^{\infty }\underset{i_{1},...,i_{n}\in \{1,...,s\}}{%
\sum }E\left( A_{i_{1}}...A_{i_{n}}\right) z_{i_{1}}...z_{i_{n}}$

$\ \ \ \ \ \ =\sum_{n=1}^{\infty }\underset{i_{1},...,i_{n}\in \{1,...,s\}}{%
\sum }\left( \varphi (P_{1}^{(i_{1},...,i_{n})}),....,\varphi
(P_{N}^{(i_{1},...,i_{n})})\right) \,z_{i_{1}}...z_{i_{n}}.$

\strut \strut

In this section, we will show that

\strut

$K_{n}\left(
(a_{1}^{(i_{1})},...,a_{N}^{(i_{1})}),...,(a_{1}^{(i_{n})},...,a_{N}^{(i_{n})})\right) 
$

\begin{center}
$=\left(
k_{n}(Q_{1}^{(i_{1},...,i_{n})}),...,k_{n}(Q_{N}^{(i_{1},...,i_{n})})\right)
,$
\end{center}

\strut

where $(k_{n})_{n=1}^{\infty }$ is the cumulants induced by $\varphi ,$ in
the sense of Nica and Speicher, and

\strut

\begin{center}
$Q_{j}^{(i_{1})}=a_{j}^{(i_{1})},$ \ $j=1,...,N,$
\end{center}

\strut

\begin{center}
$Q_{j}^{(i_{1},i_{2})}=\sum_{k=1}^{2}\left(
a_{k}^{(i_{1})},a_{(j+1)-k}^{(i_{2})}\right) \in A\times A,$ \ $j=1,...,N$
\end{center}

and

\strut

\begin{center}
$Q_{j}^{(i_{1},...,i_{n-1},i_{n})}=\sum_{k=1}^{n}\left(
Q_{j}^{(i_{1},...,i_{n-1})},a_{(j+1)-k}^{(i_{n})}\right) \in \,\underset{%
n-times}{\underbrace{A\times ....\times A}}$ ,
\end{center}

\strut

for all $j=1,...,N.$ Actually, we can regard $\underset{n-times}{\underbrace{%
A\times ...\times A}}$ \ as $\ \underset{n-times}{\underbrace{A\otimes
_{B}....\otimes _{B}A}},$ \ since $C_{T}(B)=B,$ whenever $B=\mathcal{C}^{N}$
and $T$ is a Toeplitz matricial probability space over $B.$ So,

$\strut \strut $

$k_{n}\left( \sum_{k=1}^{n}\left(
Q_{j}^{(i_{1},...,i_{n-1})},a_{(j+1)-k}^{(i_{n})}\right) \right) $

$\ \ \ \ \ \ \ \ \ =\sum_{k=1}^{n}k_{n}\left(
Q_{j}^{(i_{1},...,i_{n-1})},a_{(j+1)-k}^{(i_{n})}\right) $

$\ \ \ \ \ \ \ \ \ =\sum_{k=1}^{n}c^{(n)}\left( Q_{j}^{\otimes
(i_{1},...,i_{n-1})}\otimes a_{(j+1)-k}^{(i_{n})}\right) ,$

\strut

is well-detemined, where $c^{(n)}$ is the $n$-thcomponent of the cumulant
multiplicative function $\widehat{c}=(c^{(k)})_{k=1}^{\infty }\in I(A,\Bbb{C}%
),$ induced by $\varphi $ and where

\strut

\begin{center}
$Q_{j}^{\otimes (i_{1},...,i_{n-1})}=Q_{j}^{\otimes
(i_{1},...,i_{n-2})}\otimes a_{(j+1)-k}^{(i_{n-1})},$
\end{center}

\strut

with

\strut

\begin{center}
$Q_{j}^{\otimes (i_{1})}=a_{j}^{(i_{1})}$ \ and \ $Q_{j}^{\otimes
(i_{1},i_{2})}=\sum_{k=1}^{2}a_{k}^{(i_{1})}\otimes a_{(j+1)-k}^{(i_{2})},$
\end{center}

\strut

for all $j=1,...,N.$

\strut

\begin{definition}
Fix $s\in \Bbb{N}.$ Let

\begin{center}
$P_{j}^{(i_{1},...,i_{n})}=%
\sum_{k=1}^{j}P_{k}^{(i_{1},...,i_{n-1})}a_{(j+1)-k}^{(i_{n})}\in (A,\varphi
)$
\end{center}

with

\begin{center}
$P_{j}^{(i_{1})}=a_{j}^{(i_{1})}$ \ and \ $P_{j}^{(i_{1},i_{2})}=%
\sum_{k=1}^{j}a_{k}^{(i_{1})}a_{(j+1)-k}^{(i_{2})},$
\end{center}

\strut

in $(A,\varphi ),$ be determined inductively, for all $j=1,...,N.$ Then,
with respect to them, we will define

\strut

$Q_{j}^{(i_{1})}=a_{j}^{(i_{1})}\in A,$

$Q_{j}^{(i_{1},i_{2})}=\sum_{k=1}^{j}\left(
a_{j}^{(i_{1})},a_{(j+1)-k}^{(i_{2})}\right) \in A\times A$

\strut

and inductively on $n$,

\strut

$Q_{j}^{(i_{1},...,i_{n})}=\sum_{k=1}^{j}\left(
Q_{j}^{(i_{1},...,i_{n-1})},a_{(j+1)-k}^{(i_{n})}\right) \in \,\underset{%
n-times}{\underbrace{A\times ...\times A}}$ ,

\strut

for all $(i_{1},...,i_{n})\in \{1,...,s\}^{n},$ $n\in \Bbb{N}$, for all $%
j=1,...,N.$
\end{definition}

\strut \strut

\begin{example}
\strut At the end of the previous section, we observed an example such that

\strut

\begin{center}
$P_{1}^{(1,2,3)}=a_{1}^{(1)}a_{1}^{(2)}a_{1}^{(3)}$
\end{center}

and

\begin{center}
$%
P_{2}^{(1,2,3)}=a_{1}^{(1)}a_{1}^{(2)}a_{2}^{(3)}+a_{1}^{(1)}a_{2}^{(2)}a_{1}^{(3)}+a_{2}^{(1)}a_{1}^{(2)}a_{1}^{(3)}. 
$
\end{center}

\strut

In this case, we can decide

\strut

\begin{center}
$Q_{1}^{(1,2,3)}=\left( a_{1}^{(1)},a_{1}^{(2)},a_{1}^{(3)}\right) $
\end{center}

and

\begin{center}
$Q_{2}^{(1,2,3)}=\left( a_{1}^{(1)},a_{1}^{(2)},a_{2}^{(3)}\right) +\left(
a_{1}^{(1)},a_{2}^{(2)},a_{1}^{(3)}\right) +\left(
a_{2}^{(1)},a_{1}^{(2)},a_{1}^{(3)}\right) .$
\end{center}

\strut

Notice that the addition above is meaningless. This addition is not a
pointwise addition ! This addition is only depending on cumulant $k_{n},$
induced by $\varphi .$ And we can get that

\strut

\begin{center}
$Q_{1}^{\otimes (1,2,3)}=a_{1}^{(1)}\otimes a_{1}^{(2)}\otimes a_{1}^{(3)}$
\end{center}

and

\strut

$Q_{2}^{\otimes (1,2,3)}=\left( a_{1}^{(1)}\otimes a_{1}^{(2)}\otimes
a_{2}^{(3)}\right) +\left( a_{1}^{(1)}\otimes a_{2}^{(2)}\otimes
a_{1}^{(3)}\right) +\left( a_{2}^{(1)}\otimes a_{1}^{(2)}\otimes
a_{1}^{(3)}\right) .$
\end{example}

\strut

\begin{proposition}
Let $B=\mathcal{C}^{N}$ be a Toeplitz matricial algebra and let $(A,\varphi
) $ be a NCPSpace. Let $(T,E)$ be a Toeplitz matricial probability space
over $B$ and let $%
A_{1}=(a_{1}^{(1)},...,a_{N}^{(1)}),...,A_{s}=(a_{1}^{(s)},...,a_{N}^{(s)})%
\in (T,E)$ be $B$-valued random variables ($s\in \Bbb{N}$). Then the
R-transform of $A_{1},...,A_{s}$ is determined by

\strut

$R_{A_{1},...,A_{s}}(z_{1},...,z_{s})\overset{def}{=}%
R_{A_{1},...,A_{s}}^{t}(z_{1},...,z_{s})$

\begin{center}
$=\sum_{n=1}^{\infty }\underset{i_{1},...,i_{n}\in \{1,...,s\}}{\sum }%
K_{n}(A_{i_{1}},...,A_{i_{n}})\,z_{i_{1}}...z_{i_{n}},$
\end{center}

\strut

where

\strut

\begin{center}
$K_{n}(A_{i_{1}},...,A_{i_{n}})=\underset{\pi \in NC(n)}{\sum }\widehat{E}%
(\pi )\left( A_{i_{1}}\otimes ...\otimes A_{i_{n}}\right) \mu (\pi ,1_{n}),$
\end{center}

\strut

for all $(i_{1},...,i_{n})\in \{1,...,s\}^{n},$ $n\in \Bbb{N}.$
\end{proposition}

\strut

\begin{proof}
By the M\"{o}bius inversion, we can always define a $(i_{1},...,i_{n})$-th
cumulants from moments. We defined that, since $C_{T}(B)=B,$ moments of $%
A_{1},...,A_{s}$ by trivial moments of $A_{1},...,A_{s}.$ So, WLOG, we can
define a $(i_{1},...,i_{n})$-th cumulants of $A_{1},...,A_{s}$ by trivial
(amalgamated) cumulants of them. So, by Section 1.2, we have that

\strut

$K_{n}\left( A_{i_{1}},...,A_{i_{n}}\right) =c^{(n)}\left( A_{i_{1}}\otimes
...\otimes A_{i_{n}}\right) $

$\ \ \ \ \ \ \ \ \ \ \ \ \ \ \ \ \ \ \ \ \ \ \ \ \ =\underset{\pi \in NC(n)}{%
\sum }\widehat{E}(\pi )(A_{i_{1}}\otimes ...\otimes A_{i_{n}})\mu (\pi
,1_{n}),$

\strut

where $\widehat{E}=(E^{(k)})_{k=1}^{\infty }\in I^{m}(T,B)$ is the moment
multiplicative bimodule map induced by the $B$-functional $E:T\rightarrow B.$
Also, again, by Section 1.2, we can define

\strut

\begin{center}
$%
R_{A_{1},...,A_{s}}(z_{1},...,z_{s})=R_{A_{1},...,A_{s}}^{t}(z_{1},...,z_{s}). 
$
\end{center}
\end{proof}

\strut

\strut

\begin{definition}
Let $B=\mathcal{C}^{N}$ be a Toeplitz matricial algebra and let $(A,\varphi
) $ be a NCPSpace with its linear functional $\varphi :A\rightarrow \Bbb{C}.$
Let $(T_{N},E_{N})\equiv (T,E)$ be a Toeplitz matricial probability space
over $B$ and let $%
(a_{1}^{(1)},...,a_{N}^{(1)}),...,(a_{1}^{(s)},...,a_{N}^{(s)})$ be $B$%
-valued random variables ($s\in \Bbb{N}$). For $(i_{1},...,i_{n})\in
\{1,...,s\}^{n},$ $n\in \Bbb{N},$ define a $(i_{1},...,i_{n})$-th cumulants,
by

\strut

$K_{n}\left(
(a_{1}^{(i_{1})},...,a_{N}^{(i_{1})}),...,(a_{1}^{(i_{n})},...,a_{N}^{(i_{n})})\right) 
$

\begin{center}
$=\underset{\pi \in NC(n)}{\sum }\widehat{E}(\pi )\left(
(a_{1}^{(i_{1})},...,a_{N}^{(i_{1})})\otimes ...\otimes
(a_{1}^{(i_{n})},...,a_{N}^{(i_{n})})\right) \mu (\pi ,1_{n}).$
\end{center}
\end{definition}

\strut

\strut

\begin{lemma}
Let $(T,E)$ be a Toeplitz matricial probability space over $B=\mathcal{C}%
^{N},$ induced by $(A,\varphi )$ and let $%
A_{1}=(a_{1}^{(1)},...,a_{N}^{(1)}),$ ..., $%
A_{s}=(a_{1}^{(s)},...,a_{N}^{(s)})\in (T,E)$ be $B=\mathcal{C}^{N}$-valued
random varialbes ($s\in \Bbb{N}$). Then $(i_{1},...,i_{n})$-th moment of $%
A_{1},...,A_{s}$ has the following relation ;

\strut

\begin{center}
$E\left( A_{i_{1}}...A_{i_{n}}\right) =\left( \underset{\pi \in NC(n)}{\sum }%
k_{\pi }(Q_{1}^{(i_{1},...,i_{n})}),...,\underset{\pi \in NC(n)}{\sum }%
k_{\pi }\left( Q_{N}^{(i_{1},...,i_{n})}\right) \right) ,$
\end{center}

\strut

for all $(i_{1},...,i_{n})\in \{1,...,s\}^{n},$ $n\in \Bbb{N}$, where $%
k_{\pi }$ is a cumulants depending on $\pi \in NC(n),$ induced by $\varphi ,$
with respect to $(A,\varphi ).$
\end{lemma}

\strut

\begin{proof}
By the previous observation, we have that the $(i_{1},...,i_{n})$-th
(trivial) moment of $A_{1},...,A_{s}$ is

\strut

$E\left( A_{i_{1}}...A_{i_{n}}\right) $

$\ \ \ =E\left(
(a_{1}^{(i_{1})},...,a_{N}^{(i_{1})})...(a_{1}^{(i_{n})},...,a_{N}^{(i_{n})})\right) 
$

$\ \ \ =E\left(
(P_{1}^{(i_{1},...,i_{n})},...,P_{N}^{(i_{1},...,i_{n})})\right) $

\strut \strut

where

\begin{center}
$P_{1}^{(i_{1},...,i_{n})}=a_{1}^{(i_{1})}a_{1}^{(i_{2})}...a_{1}^{(i_{n})}$
\end{center}

and

\begin{center}
$P_{j}^{(i_{1},...,i_{n})}=%
\sum_{k=1}^{j}P_{k}^{(i_{1},...,i_{n-1})}a_{(j+1)-k}^{(i_{n})},$ for $%
j=2,...,N,$
\end{center}

hence

\strut

$\ \ \ =\left( \varphi (P_{1}^{(i_{1},...,i_{n})}),...,\varphi
(P_{N}^{(i_{1},...,i_{n})})\right) \in B.$

\strut

By the M\"{o}bius inversion,

\strut

\begin{center}
$
\begin{array}{ll}
\varphi \left( P_{1}^{(i_{1},...,i_{n})}\right) & =\varphi \left(
a_{1}^{(i_{1})}a_{1}^{(i_{2})}...a_{1}^{(i_{n})}\right) \\ 
& =\underset{\pi \in NC(n)}{\sum }k_{\pi }\left(
a_{1}^{(i_{1})},a_{1}^{(i_{2})},...,a_{1}^{(i_{n})}\right)
\end{array}
$
\end{center}

and

\begin{center}
$\varphi \left( P_{j}^{(i_{1},...,i_{n})}\right) =\underset{\theta \in NC(n)%
}{\sum }k_{\theta }\left( Q_{j}^{(i_{1},...,i_{n})}\right) ,$
\end{center}

\strut

for all $j=2,...,N,$ by the definition of $Q_{j}^{(i_{1},...,i_{n})},$ where
above $k_{\pi }$ and $k_{\theta },$ for $\pi ,\theta \in NC(n),$ are defined
in the sense of Nica and Speicher (See Section 1.3 or see [2] and [3]). Thus

\strut

$E\left( A_{i_{1}}...A_{i_{n}}\right) =\left( \varphi
(P_{1}^{(i_{1},...,i_{n})}),...,\varphi (P_{N}^{(i_{1},...,i_{n})})\right) $

$\ \ \ \ \ \ \ \ \ =\left( \underset{\pi \in NC(n)}{\sum }k_{\pi }\left(
Q_{1}^{(i_{1},...,i_{n})}\right) ,...,\underset{\pi \in NC(n)}{\sum }k_{\pi
}\left( Q_{N}^{(i_{1},...,i_{n})}\right) \right) \in B,$

\strut

for all $(i_{1},...,i_{n})\in \{1,...,s\}^{n},$ $n\in \Bbb{N}.$
\end{proof}

\strut

\strut

\begin{theorem}
Let $B=\mathcal{C}^{N}$ be a Toeplitz matricial algebra and let $(A,\varphi
) $ be a NCPSpace with its linear functional $\varphi :A\rightarrow \Bbb{C}.$
Let $(T,E)$ be a Toeplitz matricial probability space over $B$ with its $B$%
-functional $E:T\rightarrow B$ and let $%
A_{1}=(a_{1}^{(1)},...,a_{N}^{(1)}),...,A_{s}=(a_{1}^{(s)},...,a_{N}^{(s)})%
\in (T,E)$ be $B$-valued random variables ($s\in \Bbb{N}$). Then the $%
(i_{1},...,i_{n})$-th cumulant of $A_{1},...,A_{s}$ is

\strut

$K_{n}\left( A_{i_{1}},...,A_{i_{n}}\right) =\left(
k_{n}(Q_{1}^{(i_{1},...,i_{n})}),...,k_{n}(Q_{N}^{(i_{1},...,i_{n})})\right)
\in B$

\strut

for all $(i_{1},...,i_{n})\in \{1,...,s\}^{n},$ $n\in \Bbb{N}.$
\end{theorem}

\strut

\begin{proof}
By definition, we have that

\strut

\begin{center}
$K_{n}\left( A_{i_{1}},...,A_{i_{n}}\right) =\underset{\pi \in NC(n)}{\sum }%
\widehat{E}(\pi )\left( A_{i_{1}}\otimes ...\otimes A_{i_{n}}\right) \mu
(\pi ,1_{n}),$
\end{center}

\strut

for any $(i_{1},...,i_{n})\in \{1,...,s\}^{n},$ $n\in \Bbb{N}.$ Fix $%
(i_{1},...,i_{n}).$ Then the right-hand side of the previous formular can be
rewritten

\strut

\begin{center}
$\underset{\pi \in NC(n)}{\sum }\,\underset{V\in \pi (o)}{\prod }\widehat{E}%
(\pi \mid _{V})\left( A_{i_{1}}\otimes ...\otimes A_{i_{n}}\right) \mu
\left( \pi \mid _{V},1_{n}\mid _{V}\right) .$
\end{center}

\strut

Recall that if $B=\mathcal{C}^{N},$ then $C_{T}(B)=B.$ Thus, for any $%
(j_{1},...,j_{k})\in \{1,...,s\}^{k},$ $k\in \Bbb{N},$

\strut

\begin{center}
$E\left( A_{j_{1}}\otimes ...\otimes A_{j_{k}}\right) \in B=C_{T}(B).$
\end{center}

\strut

By the previous property, also we can rewrite the right-hand side of the
previous formular, by

\strut

$\underset{\pi \in NC(n)}{\sum }\,\underset{W\in \pi =\pi (o)\cup \pi (i)}{%
\prod }\widehat{E}(\pi \mid _{W})\left( A_{i_{1}}\otimes ...\otimes
A_{i_{n}}\right) \mu (\pi \mid _{W},1_{n}\mid _{W}).$------(*)

\strut

Now, let $W=(j_{1},...,j_{k})\in \pi =\pi (o)\cup \pi (i)$ be a block, where 
$\pi \in NC(n).$ Then, by the previous lemma, we have that

\strut

$\widehat{E}(\pi \mid _{W})\left( A_{i_{1}}\otimes ...\otimes
A_{i_{n}}\right) =E^{(k)}\left( A_{j_{1}}\otimes ...\otimes A_{j_{k}}\right) 
$

$\ \ \ \ =E\left( A_{j_{1}}...A_{j_{k}}\right) $

$\ \ \ \ =E\left(
(a_{1}^{(j_{1})},...,a_{N}^{(j_{1})})...(a_{1}^{(j_{k})},....,a_{N}^{(j_{k})})\right) 
$

$\ \ \ \ =E\left(
(P_{1}^{(j_{1},...,j_{k})},...,P_{N}^{(j_{1},...,j_{k})})\right) $

$\ \ \ \ =\left( \varphi (P_{1}^{(j_{1},...,j_{k})}),....,\varphi
(P_{N}^{(j_{1},...,j_{k})})\right) $

$\ \ \ \ =\left( \underset{\theta \in NC(k)}{\sum }k_{\theta
}(Q_{1}^{(j_{1},...,j_{k})}),....,\underset{\theta \in NC(k)}{\sum }%
k_{\theta }(Q_{N}^{(j_{1},...,j_{k})})\right) \in B.$

\strut \strut

Therefore,

\strut

\begin{center}
(*) $=\underset{\pi \in NC(n)}{\sum }\,\underset{W\in \pi }{\prod }\left( 
\underset{\theta \in NC(k)}{\sum }k_{\theta }(Q_{1}^{(W)}),....,\underset{%
\theta \in NC(k)}{\sum }k_{\theta }(Q_{N}^{(W)})\right) \mu (\pi \mid
_{W},1_{n}\mid _{W}),$
\end{center}

\strut

where $Q_{j}^{(W)}=Q_{j}^{(j_{1},...,j_{k})},$ if $W=(j_{1},...,j_{k})\in
\pi \in NC(n)$ \ ($k<n$). Equivalently, we have that

\strut

\begin{center}
$K_{n}\left( A_{i_{1}},...,A_{i_{n}}\right) =\left(
k_{n}(Q_{1}^{(i_{1},...,i_{n})}),....,k_{n}(Q_{N}^{(i_{1},...,i_{n})})%
\right) .$
\end{center}
\end{proof}

\strut

The proof of the above theorem is motivated by cumulants of type B studied
in [31] (See Section 6.2 and 6.3 in [31]).

\strut

\begin{example}
We will use the same example which we have observed. Suppose we are given $%
A_{1}=(a_{1}^{(1)},a_{2}^{(1)}),$ $A_{2}=(a_{1}^{(2)},a_{2}^{(2)})$ and $%
A_{3}=(a_{1}^{(3)},a_{2}^{(3)}),$ $\mathcal{C}^{2}$-valued random variables
in $(T_{2},E_{2}).$ Then we have that

\strut

\begin{center}
$Q_{1}^{(1,2,3)}=\left( a_{1}^{(1)},a_{1}^{(2)},a_{1}^{(3)}\right) $
\end{center}

and

\begin{center}
$Q_{2}^{(1,2,3)}=\left( a_{1}^{(1)},a_{1}^{(2)},a_{2}^{(3)}\right) +\left(
a_{1}^{(1)},a_{2}^{(2)},a_{1}^{(3)}\right) +\left(
a_{2}^{(1)},a_{1}^{(2)},a_{1}^{(3)}\right) .$
\end{center}

\strut

So,

\strut

$K_{3}\left( A_{1},A_{2},A_{3}\right) =\left(
k_{3}(Q_{1}^{(1,2,3)}),\,\,k_{3}(Q_{2}^{(1,2,3)})\right) $

\strut

$=(k_{3}(a_{1}^{(1)},a_{1}^{(2)},a_{1}^{(3)}),\,\,$

$\ \ \ \ \ \ \ \ k_{3}\left( a_{1}^{(1)},a_{1}^{(2)},a_{2}^{(3)}\right)
+k_{3}\left( a_{1}^{(1)},a_{2}^{(2)},a_{1}^{(3)}\right) +k_{3}\left(
a_{2}^{(1)},a_{1}^{(2)},a_{1}^{(3)}\right) ).$
\end{example}

\strut

\strut

\strut

\strut

\strut

\subsection{Toeplitz Matricial R-transform Theory}

\strut

\strut

\strut

\strut

Let $(A,\varphi )$ be a NCPSpace, in the sense of Section 1.3, with its
linear functional $\varphi :A\rightarrow \Bbb{C}.$ In this section, by using
information about cumulants, we will consider R-transform Calculus on a
Toeplitz matricial probability space $(T,E),$ over a Toeplitz matricial
algebra, $\ B=\mathcal{C}^{N},$ induced by a NCPSpace $(A,\varphi )$. Also,
we will observe the freeness of subsets of $T$. To do that, we need to
observe the freeness on $(T,E).$ Clearly, by the Speicher's
characterization, we can define freeness on $(T,E).$ We want to express this
freeness with respect to the freeness on $(A,\varphi ).$

\strut

Let $%
X_{1}=(x_{1}^{(1)},...,x_{N}^{(1)}),...,X_{s}=(x_{1}^{(s)},...,x_{N}^{(s)}),$
$Y_{1}=(y_{1}^{(1)},...,y_{N}^{(1)}),...,$ and $%
Y_{t}=(y_{1}^{(t)},...,y_{N}^{(t)})$ be $B=\mathcal{C}^{N}$-valued random
variables ($s,t\in \Bbb{N}$). Define two subsets of $T,$

\strut

\begin{center}
$\Xi =\{X_{1},...,X_{s}\}$ \ and \ $\Upsilon =\{Y_{1},...,Y_{t}\}.$
\end{center}

\strut

By Speicher's (amalgamated) Freeness characterization, we can say that

\strut

[Two subsets $\Xi $ and $\Upsilon $ are free over $B$]\strut

\begin{center}
$\Longleftrightarrow $[All mixed $B$-valued cumulants vanish].
\end{center}

Remark that, since $C_{T}(B)=B,$ ''all mixed cumulants'' can be replaced by
''all \textbf{trivial} mixed cumulants'' and these trivial cumulants is
defined by $B=\mathcal{C}^{N}$-valued cumulants, above. Also, remark that,
since $C_{T}(B)=B,$ we can avoid the insertion property when we compute the $%
B$-valued cumulants. By the previous section, if $A_{i_{1}},...,A_{i_{n}}\in
(T,E)$ are $B$-valued random variables, then we have that

\strut

$K_{n}\left( A_{i_{1}},...,A_{i_{n}}\right) \overset{def}{=}\underset{\pi
\in NC(n)}{\sum }\widehat{E}(\pi )\left( A_{i_{1}}\otimes ...\otimes
A_{i_{n}}\right) \mu (\pi ,1_{n})$

$\ =\underset{\pi \in NC(n)}{\sum }\,\underset{V\in \pi }{\prod }\widehat{E}%
(\pi \mid _{V})\left( A_{i_{1}}\otimes ...\otimes A_{i_{n}}\right) \mu (\pi
\mid _{V},1_{\left| V\right| })$

$\ =\underset{\pi \in NC(n)}{\sum }\,\underset{V\in \pi }{\prod }\left(
\varphi (P_{1}^{(i_{1},...,i_{n})}),....,\varphi
(P_{N}^{(i_{1},...,i_{n})})\right) \mu (\pi \mid _{V},1_{\left| V\right| })$

$\ =\underset{\pi \in NC(n)}{\sum }\,\underset{V\in \pi }{\prod }\left( 
\underset{\theta \in NC(\left| V\right| )}{\sum }k_{\theta
}(Q_{1}^{(V)}),....,\underset{\theta \in NC(\left| V\right| )}{\sum }%
k_{\theta }(Q_{N}^{(V)})\right) \mu (\pi \mid _{V},1_{\left| V\right| })$

$\ =\left(
k_{n}(Q_{1}^{(i_{1},...,i_{n})}),...,k_{n}(Q_{N}^{(i_{1},...,i_{n})})\right)
\in B,$

\strut

where $K_{n}$'s are $B$-valued cumulants induced by $E$ and $k_{n}$'s are
cumulants, in the sense of Nica and Speicher, induced by $\varphi $ and where

\strut

\begin{center}
$P_{j}^{(i_{1},...,i_{n})}=%
\sum_{k=1}^{j}P_{k}^{(i_{1},...,i_{n-1})}a_{(j+1)-k}^{(i_{n})}\in (A,\varphi
),$
\end{center}

and

\begin{center}
$Q_{j}^{(i_{1},...,i_{n})}=\sum_{k=1}^{j}\left(
Q_{k}^{(i_{1},...,i_{n-1})},a_{(j+1)-k}^{(i_{n})}\right) \in \,\underset{%
n-times}{\underbrace{A\times ...\times A}}$ ,
\end{center}

\strut

for all $j=1,...,N.$

\strut

\begin{theorem}
Let $B=\mathcal{C}^{N}$ be a Toeplitz matricial algebra and let $(A,\varphi
) $ be an arbitrary NCPSpace. Let $(T_{N},E_{N})\equiv (T,E)$ be a Toeplitz
matricial probability space over $B$ and let $%
X_{1}=(x_{1}^{(1)},...,x_{N}^{(1)}),...,X_{s}=(x_{1}^{(s)},...,x_{N}^{(s)}),Y_{1}=(y_{1}^{(1)},...,y_{N}^{(1)}),..., 
$ and $Y_{t}=(y_{1}^{(t)},...,y_{N}^{(t)})$ be $B$-valued random variables ($%
s,t\in \Bbb{N}$). Define subsets $\Xi =\{X_{1},...,X_{s}\}$ and $\Upsilon
=\{Y_{1},...,Y_{t}\}$ in $T.$ Then $\Xi $ and $\Upsilon $ are free over $B,$
in $(T,E)$ if and only if
\end{theorem}

\strut

\begin{proof}
Suppose that $\Xi $ and $\Upsilon $ are free over $B,$ in $(T,E).$ Then, by
definition, all $B$-valued mixed cumulants of $\Xi $ and $\Upsilon $ vanish ;

\strut

\begin{center}
$K_{n}\left( A_{i_{1}},...,A_{i_{n}}\right) =\left( k_{n}\left(
Q_{1}^{(i_{1},...,i_{n})}\right) ,...,k_{n}\left(
Q_{N}^{(i_{1},...,i_{n})}\right) \right) =0_{B},$
\end{center}

\strut

where

\begin{center}
$0_{B}\equiv \left( 
\begin{array}{lll}
0 & \cdots & 0 \\ 
\vdots & \ddots & \vdots \\ 
0 & \cdots & 0
\end{array}
\right) _{N\times N}\equiv (0,...,0)\in B=\mathcal{C}^{N}$
\end{center}

\strut

and $A_{i_{1}},...,A_{i_{n}}\in \Xi \cup \Upsilon $ (mixed) in $(T,E).$ It
suffices to consider the case when $A_{i_{j}}=Y_{i_{j}}$ and $%
A_{i_{k}}=X_{i_{k}},$ $k\neq j$ in $\{1,...,n\}.$ i.e put

\strut

\begin{center}
$A_{i_{1}}=X_{i_{1}},...,A_{i_{j-1}}=X_{i_{j-1}},$ \frame{$%
A_{i_{j}}=Y_{i_{j}}$} $,A_{i_{j+1}}=X_{i_{j+1}}...,A_{i_{n}}=X_{i_{n}}.$
\end{center}

\strut

First, we will consider $P_{1}^{(i_{1},...,i_{n})}$,..., $%
P_{N}^{(i_{1},...,i_{n})},$ because $Q_{j}^{(i_{1},...,i_{n})}$ 's are
determined by them ;

\strut

\begin{center}
$%
P_{1}^{(i_{1},...,i_{n})}=x_{1}^{(i_{1})}...x_{1}^{(i_{j-1})}y_{1}^{(i_{j})}x_{1}^{(i_{j+1})}...x_{1}^{(i_{n})} 
$
\end{center}

\strut and

\begin{center}
$P_{j}^{(i_{1},...,i_{n})}=%
\sum_{k=1}^{j}P_{j}^{(i_{1},...,i_{n-1})}a_{(j+1)-k}^{(i_{n})},$ for all $%
j=1,...,N.$
\end{center}

\strut

Note that if $i_{j}=i_{n},$ then $P_{j}^{(i_{1},...,i_{n-1},i_{n})}$ has the
form of sum of mixed words in $%
x_{1}^{(i_{1})},...,x_{N}^{(i_{1})},...,x_{1}^{(i_{n-1})},...,x_{N}^{(i_{n-1})} 
$ and \frame{$%
y_{1}^{(i_{j})}=y_{1}^{(i_{n})},...,y_{N}^{(i_{j})}=y_{N}^{(i_{n})}$}$.$ If $%
1<j<n,$ then, also, $P_{j}^{(i_{1},...,i_{j},...,i_{n})}$ has the form of
sum of mixed words \ in $%
x_{1}^{(i_{1})},...,x_{N}^{(i_{1})},...,x_{1}^{(i_{j-1})},...,x_{N}^{(i_{j-1})}, 
$ \frame{$y_{1}^{(i_{j})},...,y_{N}^{(i_{j})}$}$,$ $%
x_{1}^{(i_{j+1})},...,x_{N}^{(i_{j+1})},$ $...,$ $%
x_{1}^{(i_{n})},...,x_{N}^{(i_{n})}.$ We can conclude, by induction. So,
under our assumption, generally, for $1\leq j\leq n,$

\strut

$K_{n}\left( A_{i_{1}},...,A_{i_{n}}\right) =\left( k_{n}\left(
Q_{1}^{(i_{1},...,i_{n})}\right) ,...,k_{n}\left(
Q_{N}^{(i_{1},...,i_{n})}\right) \right) $

$\ \ \
=(k_{n}(x_{1}^{(i_{1})},...,y_{1}^{(i_{j})},...,x_{1}^{(i_{n})}),\,%
\sum_{k=1}^{2}k_{n}(Q_{2}^{(i_{1},...,i_{n-1})},a_{(2+1)-k}^{(i_{n})}),$

$\ \ \ \ \ \ \ \ \ \ \ \ \ \ \ \ \ \ \ \ \ \ \ \
...,%
\sum_{k=1}^{N}k_{n}(Q_{N}^{(i_{1},...,i_{n-1})},a_{(j+1)-k}^{(i_{n})}))=(0,...,0)\in B. 
$

\strut

$\Longleftrightarrow $

\strut

$%
k_{n}(x_{1}^{(i_{1})},...,x_{1}^{(i_{j-1})},y_{1}^{(i_{j})},x_{1}^{(i_{j+1})},...,x_{1}^{(i_{n})})=0\in 
\Bbb{C}$,

\strut

$\sum_{k=1}^{2}k_{n}(Q_{2}^{(i_{1},...,i_{n-1})},a_{(2+1)-k}^{(i_{n})})=0\in 
\Bbb{C},$

$\ \ \ \ \ \ \ \ \ \ \ \ \ \ \ \ \vdots $

$\sum_{k=1}^{N}k_{n}(Q_{N}^{(i_{1},...,i_{n-1})},a_{(j+1)-k}^{(i_{n})})=0\in 
\Bbb{C}.$

\strut

$\Longleftrightarrow $

\strut

All mixed cumulants of $%
\{x_{1}^{(i_{1})},...,x_{N}^{(i_{1})},...,x_{1}^{(i_{j-1})},...,x_{N}^{(i_{j-1})},x_{1}^{(i_{j+1})},...,x_{N}^{(i_{j+1})},...,x_{1}^{(i_{n})},...,x_{N}^{(i_{n})}\} 
$ and $\{y_{1}^{(i_{j})},...,y_{N}^{(i_{j})}\}$ are free (over $\Bbb{C}$),
in the sense of Section 1.3, in $(A,\varphi ).$ Therefore, we can get that

\strut

$\Xi
=\{X_{1}=(x_{1}^{(1)},...,x_{N}^{(1)}),...,X_{s}=(x_{1}^{(s)},...,x_{N}^{(s)})\} 
$ and $\Upsilon
=\{Y_{1}=(y_{1}^{(1)},...,y_{N}^{(1)}),...,Y_{t}=(y_{1}^{(t)},...,y_{N}^{(t)})\} 
$ are free over $B=\mathcal{C}^{N},$ in $(T,E)$ if and only if $%
\{x_{1}^{(l)},...,x_{N}^{(l)}:l=1,...,s\}$ and $\{y_{1}^{(l^{\prime
})},...,y_{N}^{(l^{\prime })}:l^{\prime }=1,...,t\}$ are free, in the sense
of Section 1.3, in $(A,\varphi ).$
\end{proof}

Let

\begin{center}
$X_{1}\equiv \left( 
\begin{array}{llll}
x_{1}^{(1)} & x_{2}^{(1)} & \cdots & x_{N}^{(1)} \\ 
& x_{1}^{(1)} & \ddots & \vdots \\ 
&  & \ddots & x_{2}^{(1)} \\ 
O &  &  & x_{1}^{(1)}
\end{array}
\right) ,...,X_{s}\equiv \left( 
\begin{array}{llll}
x_{1}^{(s)} & x_{2}^{(s)} & \cdots & x_{N}^{(s)} \\ 
& x_{1}^{(s)} & \ddots & \vdots \\ 
&  & \ddots & x_{2}^{(s)} \\ 
O &  &  & x_{1}^{(s)}
\end{array}
\right) $
\end{center}

and

\begin{center}
$Y_{1}\equiv \left( 
\begin{array}{llll}
y_{1}^{(1)} & y_{2}^{(1)} & \cdots & y_{N}^{(1)} \\ 
& y_{1}^{(1)} & \ddots & \vdots \\ 
&  & \ddots & y_{2}^{(1)} \\ 
O &  &  & y_{1}^{(1)}
\end{array}
\right) ,...,Y_{t}=\left( 
\begin{array}{llll}
y_{1}^{(t)} & y_{2}^{(t)} & \cdots & y_{N}^{(t)} \\ 
& y_{1}^{(t)} & \ddots & \vdots \\ 
&  & \ddots & y_{2}^{(t)} \\ 
O &  &  & y_{1}^{(t)}
\end{array}
\right) $
\end{center}

\strut

be $\mathcal{C}^{N}$-valued random variables in $(T_{N},E_{N}),$ induced by $%
(A,\varphi ).$ Then, by the previous theorem, we characterize $\mathcal{C}%
^{N}$-valued freeness with respect to (scarlar-valued) freeness on $%
(A,\varphi ).$ i.e $\{X_{1},...,X_{s}\}$ and $\{Y_{1},...,Y_{t}\}$ are free
over $\mathcal{C}^{N}$ if and only if $%
\{x_{1}^{(l)},...,x_{N}^{(l)}:l=1,...,s\}$ and $\{y_{1}^{(l^{\prime
})},...,y_{N}^{(l^{\prime })}:l^{\prime }=1,...,t\}$ are free in $(A,\varphi
).$

\strut

Now, we can characterize $\mathcal{C}^{N}$-valued freeness on $(T_{N},E_{N})$
with respect to scalar-valued freeness on $(A,\varphi ),$ where $%
(T_{N},E_{N})$ is a Toeplitz matricial probability space over a Toeplitz
matricial algebra, induced by $(A,\varphi ).$ So, we can observe R-transform
Calculus.

\strut

\begin{proposition}
Let $(A,\varphi )$ be a NCPSpace and let $%
x_{1}^{(l)},...,x_{N}^{(l)},y_{1}^{(l^{\prime })},...,y_{N}^{(l^{\prime
})}\in (A,\varphi )$ be random variables, where $l=1,...,s$ and $l^{\prime
}=1,...,t$ ($s,t\in \Bbb{N}$). Let $B=\mathcal{C}^{N}$ be a Toeplitz
matricial algebra and let $(T,E)$ be a NCPSpace over $B,$ induced by $%
(A,\varphi ).$ Let

\strut

\begin{center}
$X_{l}=(x_{1}^{(l)},...,x_{N}^{(l)})$ \ and \ $Y_{l^{\prime
}}=(y_{1}^{(l^{\prime })},...,y_{N}^{(l^{\prime })})$
\end{center}

\strut

be $B$-valued random variables in $(T,E),$ where $l=1,...,s$ and $l^{\prime
}=1,...,t.$ If $\widetilde{X}=\{x_{1}^{(l)},...,x_{N}^{(l)}:l=1,...,s\}$ and 
$\widetilde{Y}=\{y_{1}^{(l^{\prime })},...,y_{N}^{(l^{\prime })}:l^{\prime
}=1,...,t\}$ are free in $(A,\varphi ),$ then

\strut

(1) $%
R_{X_{1},...,X_{s},Y_{1},...,Y_{t}}(z_{1},...,z_{s+t})=R_{X_{1},...,X_{s}}(z_{1},...,z_{s})+R_{Y_{1},...,Y_{t}}(z_{s+1},...,z_{s+t}) 
$

\strut

(2) if $s=t,$ then

\strut

$\ \ \ \
R_{X_{1}+Y_{1},...,X_{s}+Y_{s}}(z_{1},...,z_{s})=R_{X_{1},...,X_{s}}(z_{1},...,z_{s})+R_{Y_{1},...,Y_{s}}(z_{1},...,z_{s}) 
$

\strut

(3) if $s=t,$ then

\strut

$\ \ \ \ R_{X_{1}Y_{1},...,X_{s}Y_{s}}(z_{1},...,z_{s})=\left(
R_{X_{1},...,X_{s}}\,\,\frame{*}\,\,R_{Y_{1},...,Y_{s}}\right)
(z_{1},...,z_{s}).$
\end{proposition}

\strut

\begin{proof}
By the previous theorem, if $\widetilde{X}$ and $\widetilde{Y}$ are free in $%
(A,\varphi ),$ then $\{X_{1},...,X_{s}\}$ and $\{Y_{1},...,Y_{t}\}$ are free
over $B,$ in $(T,E).$ So, (1) and (2) are easily verified. We will only
observe (3), when $s=t.$ Fix $(i_{1},...,i_{n})\in \{1,...,s\}^{n},$ $n\in 
\Bbb{N}.$ Then

\strut

\begin{center}
$K_{n}\left( X_{i_{1}}Y_{i_{1}},...,X_{i_{n}}Y_{i_{n}}\right) =\left(
k_{n}(Q_{1}^{(i_{1},...,i_{n})}),...,k_{n}(Q_{N}^{(i_{1},...,i_{n})})\right)
,$
\end{center}

\strut

where

\strut

\begin{center}
$
\begin{array}{ll}
X_{i_{k}}Y_{i_{k}} & =(x_{1}^{(i_{k})},...,x_{N}^{(i_{k})})\cdot
(y_{1}^{(i_{k})},...,y_{N}^{(i_{k})}) \\ 
& =\left(
x_{1}^{(i_{k})}y_{1}^{(i_{k})},\,%
\sum_{k=1}^{2}x_{k}^{(i_{1})}y_{(2+1)-k}^{(i_{k})},...,%
\sum_{k=1}^{N}x_{k}^{(i_{n})}x_{(N+1)-k}^{(i_{n})}\right) \\ 
& \equiv (a_{1}^{(i_{k})},...,a_{N}^{(i_{k})}),
\end{array}
$
\end{center}

\strut

for all $k=1,...,n$ and

\strut

\begin{center}
$Q_{j}^{(i_{1},...,i_{n})}=\sum_{k=1}^{j}\left(
Q_{j}^{(i_{1},...,i_{n-1})},a_{(j+1)-k}^{(i_{n})}\right) \in \,\underset{%
n-times}{\underbrace{A\times ...\times A}},$
\end{center}

\strut

for all $j=1,...,N.$ Thus,

\strut

$k_{n}\left( Q_{j}^{(i_{1},...,i_{n})}\right) =\sum_{k=1}^{j}k_{n}\left(
Q_{j}^{(i_{1},...,i_{n-1})},a_{(j+1)-k}^{(i_{n})}\right) $

\strut

and, inductively, $\left(
Q_{j}^{(i_{1},...,i_{n-1})},a_{(j+1)-k}^{(i_{n})}\right) $ is the sum of
alternating tuples of $x_{p}^{(q)}$'s and $y_{p^{\prime }}^{(q^{\prime })}$%
's. By little abuse of notation, we have that

\strut

$\left(
k_{n}(Q_{1}^{(i_{1},...,i_{n})}),...,k_{n}(Q_{N}^{(i_{1},...,i_{n})})\right) 
$

$\ \ =(\underset{\theta \in NC(n)}{\sum }k_{\theta
}(Q_{1}(x_{p}^{(q)}))\cdot k_{Kr(\theta )}(Q_{1}(y_{p^{\prime }}^{(q^{\prime
})})),...$

\begin{center}
$...,\underset{\theta \in NC(n)}{\sum }k_{\theta }(Q_{N}(x_{p}^{(q)})\cdot
k_{Kr(\theta )}(Q_{N}(y_{p^{\prime }}^{(q^{\prime })}))),$
\end{center}

\strut

where $\theta \cup Kr(\theta )\in NC(2n)$ is the alternating union of
noncrossing partitions such that $x_{p}^{(q)}$'s are depending on $\theta $
and $y_{p^{\prime }}^{(q^{\prime })}$'s are depending on $Kr(\theta ).$
Therefore,

\strut

\begin{center}
$K_{n}\left( X_{i_{1}}Y_{i_{1}},...,X_{i_{n}}Y_{i_{n}}\right) =\underset{\pi
\in NC(n)}{\sum }\left( K_{\pi }(X_{i_{1}},...,X_{i_{n}})\right) \left(
K_{Kr(\pi )}(Y_{i_{1}},...,Y_{i_{n}})\right) .$
\end{center}

\strut

(Again, remark that $C_{T}(B)=B$ !!)
\end{proof}

\strut

\strut

\strut

\strut

\strut

\strut

\section{Applications}

\strut \strut

\strut

\strut

\strut

\subsection{Operator-Valued Moments and Cumulants of a Single Toeplitz
matricial valued random variable}

\strut

\strut

\strut

\strut

Let $B=\mathcal{C}^{N}$ be a Toeplitz matricial algebra and let $(A,\varphi
) $ be a NCPSpace, in the sense of Section 1.3. Let $(T,E)$ be a Toeplitz
matricial probability space over $B,$ induced by $(A,\varphi ).$ In this
section, we will compute operator-valued moments and operator-valued
cumulants of a single $B$-valued random variable

\begin{center}
$A=(a_{1},...,a_{N})\equiv \left( 
\begin{array}{llll}
a_{1} & a_{2} & \cdots & a_{N} \\ 
& a_{1} & \ddots & \vdots \\ 
&  & \ddots & a_{2} \\ 
O &  &  & a_{1}
\end{array}
\right) \in (T,E).$
\end{center}

\strut \strut

We will compute the moment series of $A,$ $M_{A}(z)$ and the R-transform of $%
A,$ $R_{A}(z),$ in $\Theta _{B}^{1}.$ To do that, we have to compute the $n$%
-th moment of $A$ and $n$-th cumulant of $A,$ as $n$-th coefficient of $%
M_{A} $ and $R_{A},$ respectively, for $n\in \Bbb{N}.$ Recall that, since $%
C_{T}(B)=B,$ we can define $M_{A}$ and $R_{A}$ by trivial moment series of $%
A $ and trivial R-transform of $A,$ in the sense of Section 1.2 (See Section
2.2 and 2.3). i.e

\strut

\begin{center}
$M_{A}(z)=\sum_{n=1}^{\infty }E\left( A^{n}\right) \,z^{n}$
\end{center}

and

\begin{center}
$R_{A}(z)=\sum_{n=1}^{\infty }K_{n}\left( A,...,A\right) \,\,z^{n},$
\end{center}

\strut

where $E\left( (a_{1},...,a_{N})\right) =\left( \varphi (a_{1}),...,\varphi
(a_{n})\right) ,$ for all $(a_{1},...,a_{N})\in (T,E)$ and $K_{n}$ is
induced by $E,$ in the sense of Section 2.3.

\strut

Consider the general expression of $A^{n},$ for $n\in \Bbb{N}.$ Similar to
the previous section, we can compute that ;

\strut

$A^{2}=\left( a_{1},...,a_{N}\right) \cdot \left( a_{1},...,a_{N}\right) $

\strut

$\ \ \ \ \equiv \left( 
\begin{array}{llll}
a_{1} & a_{2} & \cdots & a_{N} \\ 
& a_{1} & \ddots & \vdots \\ 
&  & \ddots & a_{2} \\ 
O &  &  & a_{1}
\end{array}
\right) \left( 
\begin{array}{llll}
a_{1} & a_{2} & \cdots & a_{N} \\ 
& a_{1} & \ddots & \vdots \\ 
&  & \ddots & a_{2} \\ 
O &  &  & a_{1}
\end{array}
\right) $

\strut

$\ \ \ \ =\left( 
\begin{array}{llll}
a_{1}^{2} & \sum_{k=1}^{2}a_{k}a_{(2+1)-k} & \cdots & 
\sum_{k=1}^{N}a_{k}a_{(N+1)-k} \\ 
& \text{ \ \ \ \ \ \ \ \ }a_{1}^{2} & \ddots & \text{ \ \ \ \ \ \ \ \ \ \ }%
\vdots \\ 
&  & \ddots & \sum_{k=1}^{2}a_{k}a_{(2+1)-k} \\ 
&  &  & \text{ \ \ \ \ \ \ \ \ }a_{1}^{2}
\end{array}
\right) $

\strut

$\ \ \ \equiv \left(
a_{1}^{2},\,\,\sum_{k=1}^{2}a_{k}a_{(2+1)-k},\,\,...,\,\,%
\sum_{k=1}^{N}a_{k}a_{(N+1)-k}\right) $

$\ \overset{denote}{=}\left( P_{1}^{(2)},P_{2}^{(2)},...,P_{N}^{(2)}\right)
. $

\strut

Inductively, we have that

\strut

\begin{center}
$P_{j}^{(1)}=a_{j},$
\end{center}

\strut and

\begin{center}
$P_{j}^{(n)}=\sum_{k=1}^{j}P_{k}^{(n-1)}a_{(j+1)-k},$
\end{center}

\strut

for all $j=1,...,N.$

\strut

Thus we can compute the $n$-th moment of $A$ ;

\strut

\begin{center}
$E\left( A^{n}\right) =E\left( (P_{1}^{(n)},...,P_{N}^{(n)})\right) =\left(
\varphi (P_{1}^{(n)}),...,\varphi (P_{N}^{(n)})\right) \in B.$
\end{center}

\strut For example, if $N=2,$ then we have that

\strut

$P_{1}^{(1)}=a_{1}$ and $P_{2}^{(1)}=a_{2}$

\strut

$P_{1}^{(2)}=a_{1}^{2}$ \ and \ $P_{2}^{(2)}=a_{1}a_{2}+a_{2}a_{1}$

\strut

$P_{1}^{(3)}=a_{1}^{3}$ \ and \ $%
P_{2}^{(3)}=a_{1}^{2}a_{2}+a_{1}a_{2}a_{1}+a_{2}a_{1}^{2}=%
\sum_{k=1}^{2}P_{k}^{(2)}a_{(2+1)-k}$

\strut

etc, and hence

\strut

\begin{center}
$
\begin{array}{ll}
E\left( A^{3}\right) & =\left( \varphi (a_{1}^{3}),\varphi \left(
a_{1}^{2}a_{2}+a_{1}a_{2}a_{1}+a_{2}a_{1}^{2}\right) \right) \\ 
& =\left( \varphi (a_{1}^{3}),\,\,\varphi (a_{1}^{2}a_{2})+\varphi
(a_{1}a_{2}a_{1})+\varphi (a_{2}a_{1}^{2})\right)
\end{array}
$
\end{center}

\strut and

\strut $K_{3}\left( A,A,A\right) =\left(
k_{3}(a_{1},a_{1},a_{1}),\,%
\,k_{3}(a_{1},a_{1},a_{2})+k_{3}(a_{1},a_{2},a_{1})+k_{3}(a_{2},a_{1},a_{1})%
\right) .$

\strut

\begin{proposition}
Let $B=\mathcal{C}^{N}$ be a Toeplitz matricial algebra and let $(A,\varphi
) $ be a NCPSpace. Let $(T,E)$ be a Toeplitz matricial probability space
over $B$ and let $A=(a_{1},...,a_{N})\in (T,E)$ be a $B$-valued random
variable. Then

\strut

\begin{center}
$M_{A}(z)=\sum_{n=1}^{\infty }\left( \varphi
(P_{1}^{(n)},...,P_{N}^{(n)})\right) \,z^{n}$
\end{center}

and

\begin{center}
$R_{A}(z)=\sum_{n=1}^{\infty }\left(
k_{n}(Q_{1}^{(n)}),...,k_{n}(Q_{N}^{n})\right) \,z^{n},$
\end{center}

\strut

in $\Theta _{B}^{1},$ where

\strut

$Q_{j}^{(n)}=\sum_{k=1}^{j}\left( Q_{k}^{(n-1)},\,a_{(j+1)-k}\right) \in \,%
\underset{n-times}{\underbrace{A\times ....\times A}},$

\strut

with $Q_{j}^{(1)}=a_{j},$ for all $j=1,...,N.$ \ \ $\square $
\end{proposition}

\strut

\strut

\begin{theorem}
Let $B=\mathcal{C}^{N}$ be a Toeplitz matricial algebra and let $(A,\varphi
) $ be a NCPSpace in the sense of Section 1.3. Let $(T,E)$ be a Teoplitz
matricial probability space over $B,$ induced by $(A,\varphi )$ and let $%
A=\left( a_{1},...,a_{N}\right) \in (T,E)$ be a $B$-valued random variable,
with $a_{1},...,a_{N}\in (A,\varphi )$ are random variables. If $%
\{a_{1}\},...,\{a_{N}\}$ are free in $(A,\varphi )$ (i.e $%
\{a_{1},...,a_{N}\} $ is a free family in $(A,\varphi )$), then

\strut

\begin{center}
$R_{A}(z)=\left( \varphi (a_{1}),...,\varphi (a_{N})\right)
z+b_{2}z^{2}+\sum_{n=3}^{\infty }\left( k_{n}(\,\underset{n-times}{%
\underbrace{a_{1},......,a_{1}}}\,),....,0\right) \,z^{n},$
\end{center}

\strut

in $\Theta _{B}^{1},$ where $b_{2}\in B$ such that

\strut

\begin{center}
$b_{2}=\left\{ 
\begin{array}{lll}
\left( k_{2}(a_{1},a_{1}),0,k_{2}(a_{2},a_{2}),0,...,0,k_{2}(a_{\frac{N}{2}%
},\,a_{\frac{N}{2}}),0\right) &  & N\text{ is even} \\ 
&  &  \\ 
\left( k_{2}(a_{1},a_{1}),0,k_{2}(a_{2},a_{2}),0,...,0,k_{2}(a_{\frac{N-1}{2}%
},\,a_{\frac{N-1}{2}})\right) &  & N\text{ is odd.}
\end{array}
\right. $
\end{center}
\end{theorem}

\strut

\begin{proof}
By the previous proposition, we have $n$-th coefficient of $R_{A}(z),$

\strut

\begin{center}
$K_{n}\left( \underset{n-times}{\underbrace{A,......,A}}\right) =\left(
k_{n}(Q_{1}^{(n)}),....,k_{n}(Q_{N}^{(n)})\right) .$
\end{center}

\strut

We need to consider $Q_{j}^{(n)},$ for $j=1,...,N.$ Since $Q_{j}^{(n)}$ is
gotten from $P_{j}^{(n)},$ $j=1,...,N,$ we will observe $P_{j}^{(n)}$ ;

\strut

\begin{center}
$P_{j}^{(1)}=a_{j},$
\end{center}

\strut

\begin{center}
$P_{j}^{(2)}=\sum_{k=1}^{j}P_{k}^{(1)}a_{(j+1)-k}$
\end{center}

and

\begin{center}
$P_{j}^{(n)}=\sum_{k=1}^{j}P_{k}^{(n-1)}a_{(j+1)-k},$
\end{center}

\strut

for all $j=1,...,N.$ Hence

\strut

\begin{center}
$Q_{j}^{(1)}=P_{j}^{(1)}=a_{j},$
\end{center}

\strut

\begin{center}
$Q_{j}^{(2)}=\sum_{k=1}^{j}\left( Q_{k}^{(1)},\,\,a_{(j+1)-k}\right) $
\end{center}

and

\begin{center}
$Q_{j}^{(n)}=\sum_{k=1}^{j}\left( Q_{k}^{(n-1)},a_{(j+1)-k}\right) ,$
\end{center}

\strut

for all $j=1,...,N.$ \strut Thus

\strut

$K_{1}(A)=\left( k_{1}(a_{1}),...,k_{1}(a_{N})\right) =\left( \varphi
(a_{1}),...,\varphi (a_{N})\right) ,$

\strut

$K_{2}(A,A)=(k_{2}(a_{1},a_{1}),\,k_{2}\left( a_{1},a_{2}\right)
+k_{2}(a_{2},a_{1}),$

\begin{center}
$...,k_{2}(a_{1},a_{N})+k_{2}(a_{2},a_{N-1})+...+k_{2}(a_{N},a_{1}))$
\end{center}

\strut

and

\strut

$K_{n}(A,...,A)$

$\ =\left(
k_{n}(a_{1},...,a_{1}),\sum_{k=1}^{2}k_{2}(Q_{k}^{(n-1)},a_{3-k}),...,%
\sum_{k=1}^{N}k_{n}(Q_{k}^{(n-1)},a_{(N+1)-k})\right) .$

\strut

(i) \ \ If $n=1,$ then we have that

\strut

\begin{center}
$K_{1}(A)=\left( k_{1}(a_{1}),k_{1}(a_{2}),...,k_{1}(a_{N})\right) =\left(
\varphi (a_{1}),...,\varphi (a_{N})\right) .$
\end{center}

\strut

(ii) \ \ If $n\geq 3,$ then, since $\left( Q_{k}^{(n-1)},a_{(j+1)-k}\right) $
is a mixed tuple over $\{a_{1},...,a_{N}\},$ for all $j=1,...,N$ and $k\in \{%
\frame{$2$},...,N\},$ we have that

\strut

$K_{n}\left( \underset{n-times}{\underbrace{A,.......,A}}\right) =\left(
k_{n}\left( \underset{n-times}{\underbrace{a_{1},......,a_{1}}}\right)
,0,....,0\right) $

\strut

$\ \ \equiv \left( 
\begin{array}{llll}
k_{n}(a_{1},...,a_{1}) &  &  & O \\ 
& k_{1}(a_{1},...,a_{1}) &  &  \\ 
&  & \ddots &  \\ 
O &  &  & k_{1}(a_{1},...,a_{1})
\end{array}
\right) \in B.$

\strut

(iii) If $n=2,$ then we have following two cases ;

\strut

\ \ \ \ \TEXTsymbol{<}Case 1\TEXTsymbol{>} Suppose that $N=2k,$ for $k\in 
\Bbb{N}.$ Then

\strut

$\ \ \ \ K_{2}(A,A)=(k_{2}(a_{1},a_{1}),\,k_{2}\left( a_{1},a_{2}\right)
+k_{2}(a_{2},a_{1}),$

$\ \ \ \ \ \ \ \ \ \ \ \ \ \ \ \ \ \ \ \ \ \ \
...,k_{2}(a_{1},a_{N})+k_{2}(a_{2},a_{N-1})+...+k_{2}(a_{N},a_{1}))$

\strut

$\ \ \ \ \ \ \ \ \ \ =\left(
k_{2}(a_{1},a_{1}),0,k_{2}(a_{2},a_{2}),0,k_{2}(a_{3},a_{3}),0,...,0,k_{2}(a_{%
\frac{N}{2}},\,a_{\frac{N}{2}}),0\right) $

\strut

$\ \ \ \ \ \ \ \ \ \ \equiv \left( 
\begin{array}{lllllll}
k_{2}(a_{1},a_{1}) & \text{ \ \ \ \ \ \ }0 &  &  & \cdots & k_{2}(a_{\frac{%
N-1}{2}},\,a_{\frac{N-1}{2}}) & \text{ \ \ \ \ }0 \\ 
& k_{2}(a_{1},a_{1}) & 0\text{ } &  &  & \text{ \ \ \ \ \ }* & k_{2}(a_{%
\frac{N-1}{2}},\,a_{\frac{N-1}{2}}) \\ 
&  &  &  & \ddots &  & \vdots \\ 
&  &  & \ddots & \ddots &  & \text{ \ \ \ \ } \\ 
&  &  &  &  &  &  \\ 
&  &  &  &  &  & \text{ \ \ \ \ }0 \\ 
O &  &  &  &  &  & k_{2}(a_{1},a_{1})
\end{array}
\right) _{N\times N}$

\strut

\ \ \ \ \TEXTsymbol{<}Case 2\TEXTsymbol{>} Suppose that $N=2k-1,$ for $k\in 
\Bbb{N}.$ Then

\strut

$\ \ \ \ \ K_{2}(A,A)=\left(
k_{2}(a_{1},a_{1}),0,k_{2}(a_{2},a_{2}),0,...,0,k_{2}(a_{\frac{N-1}{2}},\,a_{%
\frac{N-1}{2}})\right) $

\strut

$\ \ \ \ \ \equiv \left( 
\begin{array}{lllllll}
k_{2}(a_{1},a_{1}) & \text{ \ \ \ \ \ \ }0 &  &  & \cdots & \text{ }0 & 
k_{2}(a_{\frac{N-1}{2}},\,a_{\frac{N-1}{2}}) \\ 
& k_{2}(a_{1},a_{1}) &  &  &  & \text{ *\ \ \ \ } & \text{ \ \ \ \ \ \ }0 \\ 
&  &  &  &  & \ddots & \vdots \\ 
&  &  & \ddots & \ddots &  &  \\ 
&  &  &  &  &  & \text{ \ \ } \\ 
&  &  &  &  &  & \text{ \ \ \ \ \ \ \ \ }0 \\ 
O &  &  &  &  &  & \text{ \ \ }k_{2}(a_{1},a_{1})
\end{array}
\right) _{N\times N}$

\strut

By (i), (ii) and (iii), we have that

\strut

$R_{A}(z)=\left( \varphi (a_{1}),...,\varphi (a_{N})\right) z$

$\ \ \ \ \ \ \ \ \ \ \ \ \ \ \ +\left(
k_{2}(a_{1},a_{1}),0,k_{2}(a_{2},a_{2}),0,...,0,k_{2}(a_{\frac{N}{2}},\,a_{%
\frac{N}{2}}),0\right) z^{2}$

$\ \ \ \ \ \ \ \ \ \ \ \ \ \ \ +\left(
k_{3}(a_{1},a_{1},a_{1}),0,...,0\right) z^{3}+...$

\strut

$\ =\left( \varphi (a_{1}),...,\varphi (a_{N})\right) z+\left(
k_{2}(a_{1},a_{1}),0,k_{2}(a_{2},a_{2}),0,...,0,k_{2}(a_{\frac{N}{2}},\,a_{%
\frac{N}{2}}),0\right) z^{2}$

$\ \ \ \ \ \ \ \ \ \ \ \ \ \ \ \ \ \ \ \ \ \ \ \ \ \ \ \ \ \ \ \ \
+\sum_{n=3}^{\infty }\left( k_{n}\left( \underset{n-times}{\underbrace{%
a_{1},......,a_{1}}}\right) ,0,....,0\right) \,z^{n},$

if $N$ is even, and

\strut

$R_{A}(z)=\left( \varphi (a_{1}),...,\varphi (a_{N})\right) z$

$\ \ \ \ \ \ \ \ \ \ \ \ \ \ \ +\left(
k_{2}(a_{1},a_{1}),0,k_{2}(a_{2},a_{2}),0,...,0,k_{2}(a_{\frac{N-1}{2}},\,a_{%
\frac{N-1}{2}})\right) z^{2}$

$\ \ \ \ \ \ \ \ \ \ \ \ \ \ \ +\left(
k_{3}(a_{1},a_{1},a_{1}),0,...,0\right) z^{3}+...$

\strut

$\ =\left( \varphi (a_{1}),...,\varphi (a_{N})\right) z+\left(
k_{2}(a_{1},a_{1}),0,k_{2}(a_{2},a_{2}),0,...,0,k_{2}(a_{\frac{N-1}{2}},\,a_{%
\frac{N-1}{2}})\right) z^{2}$

$\ \ \ \ \ \ \ \ \ \ \ \ \ \ \ \ \ \ \ \ \ \ \ \ \ \ \ \ \ \ \ \ \
+\sum_{n=3}^{\infty }\left( k_{n}\left( \underset{n-times}{\underbrace{%
a_{1},......,a_{1}}}\right) ,0,....,0\right) \,z^{n},$

\strut

if $N$ is odd.
\end{proof}

\strut

\strut

\strut

\strut

\strut

\subsection{\strut $\mathcal{C}^{N}$-Evenness}

\strut

\strut

\strut

\strut

\strut

First, we will consider the $B$-valued evenness, where $B$ is an arbitrary
unital algebra. Let $B$ be a unital algebra and let $(A,E)$ be a NCPSpace
over $B.$ In [12], we defined $B$-valued evenness and we observed some
properties of it. In this section, we want to characterize the $\mathcal{C}%
^{N}$-evenness, where $B=\mathcal{C}^{N},$ a Toeplitz matricial algebra, for 
$N\in \Bbb{N}.$ Throughout this section, fix $N\in \Bbb{N}.$

\strut \strut \strut

\begin{definition}
Let $B$ be a unital algebra and let $(A,E)$ be a NCPSpace over $B.$ We say a 
$B$-valued random variable, $x\in (A,E),$ is an ($B$-valued) even element if 
$x$ satisfies the following cumulant relation ;

\strut

\begin{center}
$K_{n}\left( \underset{n-times}{\underbrace{x,......,x}}\right) =0_{B},$
whenever $n\in \Bbb{N}$ is odd,
\end{center}

\strut

where

\strut

\begin{center}
$K_{n}(x,...,x)=c^{(n)}\left( x\otimes b_{i_{2}}x\otimes ...\otimes
b_{i_{n}}x\right) ,$
\end{center}

\strut

is a $n$-th cumulant induced by a $B$-functional, $E:A\rightarrow B,$ for $%
b_{i_{2}},...,b_{i_{n}}\in B,$ arbitrary.
\end{definition}

\strut

In [12], we used the alternating definition for evenness as follows ;

\strut

\begin{center}
$x\in (A,E)$ is even $\Leftrightarrow $ moments $E\left(
xb_{2}x...b_{n}x\right) =0_{B},$ whenever $n$ is odd.
\end{center}

Also, in [12], we observed that our definition with respect to cumulnats and
the above definition with respect to moments are equivalent.

\strut

\begin{theorem}
(Also see [12])\strut Let $B$ be a unital algebra and let $(A,E)$ be a
NCPSpace over $B.$ If a $B$-valued random variable, $x\in (A,E),$ is even,
then

\strut

\begin{center}
$R_{x}(z)=\sum_{n=1}^{\infty }K_{2n}\left( \underset{2n-times}{\underbrace{%
x,......,x}}\right) \,z^{n}\in \Theta _{B}^{1}$
\end{center}

with

\strut

\begin{center}
$K_{2n}\left( \underset{2n-times}{\underbrace{x,......,x}}\right) =\underset{%
\pi \in NC^{(even)}(2n)}{\sum }\,\widehat{E}(\pi )\left( x\otimes
b_{2}x\otimes ...\otimes b_{2n}x\right) \,\mu (\pi ,1_{2n}),$
\end{center}

\strut

where

\strut

\begin{center}
$NC^{(even)}(2n)\overset{def}{=}\{\theta \in NC(2n):V\in \theta
\Leftrightarrow \left| V\right| $ is even$\}.$
\end{center}
\end{theorem}

\strut

\begin{proof}
Remark that $x$ is even if and only if that all odd moments vanish. Suppose
that $\pi \in NC(2n)$ and assume that there exists at least one $V\in \pi $
satisfies that $\left| V\right| \in \Bbb{N}$ is an odd number (i.e a block, $%
V\in \pi ,$ is an odd block).

\strut

(i) Suppose that such odd block $V$ is contained in $\pi (o).$ Then

\strut

$\widehat{E}(\pi )\left( x\otimes b_{2}x\otimes ...\otimes b_{n}x\right) $

$\ \ \ \ =\underset{W\in \pi (o)}{\prod }\widehat{E}(\pi \mid _{W})(x\otimes
b_{2}x\otimes ...\otimes b_{n}x)$

$\ \ \ \ =\left( \widehat{E}(\pi \mid _{V})(x\otimes b_{2}x\otimes
...\otimes b_{n}x)\right) $

\begin{center}
$\cdot \left( \underset{W\neq V\in \pi (o)}{\prod }\widehat{E}(\pi \mid
_{W})(x\otimes b_{2}x\otimes ...\otimes b_{n}x)\right) $
\end{center}

$\ \ \ \ =0_{B}.$

\strut

(ii) Suppose that such odd block $V$ is contained in $\pi (i).$ Then

\strut

$\widehat{E}(\pi )\left( x\otimes b_{2}x\otimes ...\otimes b_{n}x\right) $

$\ \ \ \ =\underset{W\in \pi (o)}{\prod }\widehat{E}(\pi \mid _{W})(x\otimes
b_{2}x\otimes ...\otimes b_{n}x).$

\strut

Assume that the odd block $V=(v_{1},...,v_{k})\in \pi (i)$ is an inner block
with its outer block $W_{0}=(w_{1},...,w_{p})\in \pi (o)$ and assume that
there exists $j\in \{1,...,p\}$ such that $w_{j}<v_{q}<w_{j+1},$ for all $%
q=1,...,k.$ Then the right hand side of the above equation is

\strut

$\ \ =\left( \widehat{E}(\pi \mid _{W_{0}})(x\otimes b_{2}x\otimes
...\otimes b_{n}x)\right) $

\begin{center}
$\cdot \left( \underset{W\neq W_{0}\in \pi (o)}{\prod }\widehat{E}(\pi \mid
_{W})(x\otimes b_{2}x\otimes ...\otimes b_{n}x)\right) $
\end{center}

$\ \ =\left( E^{(p)}(b_{w_{1}}x\otimes ...\otimes b_{j}x\otimes
E^{(k)}(b_{v_{1}}x\otimes ...\otimes b_{v_{k}}x)b_{j+1}x\otimes ...\otimes
b_{p}x)\right) $

\begin{center}
$\cdot \left( \underset{W\neq W_{0}\in \pi (o)}{\prod }\widehat{E}(\pi \mid
_{W})(x\otimes b_{2}x\otimes ...\otimes b_{n}x)\right) $
\end{center}

$\ \ =0_{B}.$

\strut

So, by (i) and (ii), if $\pi \in NC(2n)$ contains an odd block, then $%
\widehat{E}(\pi )(x\otimes b_{2}x\otimes ...\otimes b_{2n}x)$ vanishs. So, $%
K_{2n}(x,...,x)$ is determind by

\strut

$K_{2n}(x,...,x)=\underset{\pi \in NC(2n)}{\sum }\widehat{E}(\pi )\left(
x\otimes b_{2}x\otimes ...\otimes b_{n}x\right) \mu (\pi ,1_{2n})$

\strut

$\ \ =\underset{\pi \in NC^{(even)}(2n)}{\sum }\widehat{E}(\pi )\left(
x\otimes b_{2}x\otimes ...\otimes b_{n}x\right) \mu (\pi ,1_{2n})$

\begin{center}
$+\underset{\pi \in NC(2n)\,\,\setminus \,NC^{(even)}(2n)}{\sum }\widehat{E}%
(\pi )\left( x\otimes b_{2}x\otimes ...\otimes b_{n}x\right) \mu (\pi
,1_{2n})$
\end{center}

\strut

$\ \ =\underset{\pi \in NC^{(even)}(2n)}{\sum }\widehat{E}(\pi )\left(
x\otimes b_{2}x\otimes ...\otimes b_{n}x\right) \mu (\pi ,1_{2n})+0_{B}.$
\end{proof}

\strut

From now, let $B=\mathcal{C}^{N}$ be a Toeplitz matricial algebra. We will
characterize the Toeplitz matricial evenness. Recall that, since $%
C_{T}(B)=B, $ we just defined Toeplitz matricial moments and Toeplitz
matricial cumulants as trivial moments and trivial cumulants.

\strut \strut

\begin{theorem}
Let $B=\mathcal{C}^{N}$ be a Toeplitz matricial algebra and let $(A,\varphi
) $ be a NCPSpace in the sense of Section 1.3. Let $(T,E)$ be a Toeplitz
matricial probability space over $B,$ induced by $(A,\varphi ).$ A $B$%
-valued random variable $X=\left( x_{1},...,x_{N}\right) \in (T,E)$ is ($B$%
-valued)\ even if and only if $E\left( \underset{n-times}{\underbrace{%
X,.....,X}}\right) =0_{B},$ whenever $n$ is odd. Equivalently, $%
X=(x_{1},...,x_{N})\in (T,E)$ is even if and only if

\strut

\begin{center}
$\varphi (x_{1})=...=\varphi (x_{N})=0,$ in $\Bbb{C}$
\end{center}

and

\begin{center}
$\sum_{k=1}^{j}\varphi \left( P_{k}^{(2l)}x_{(j+1)-k}\right) =0,$ in $\Bbb{C}%
,$
\end{center}

for all $j=1,...,N,$ where

\strut

\begin{center}
$X^{2l+1}=(x_{1},...,x_{N})^{2l+1}=\left(
P_{1}^{(2l+1)},...,P_{N}^{(2l+1)}\right) ,$ $\forall \,\,l\in \Bbb{N}.$
\end{center}
\end{theorem}

\strut

\begin{proof}
By an alternating definition of amalgamated evenness, we have that $X=\left(
x_{1},...,x_{N}\right) \in (T,E)$ is even if and only if all odd moments
vanishs. So, $E\left( X^{n}\right) =0_{B},$ whenever $n\in \Bbb{N}$ is odd.
Suppose that $n$ is odd. Then

\strut

\begin{center}
$
\begin{array}{ll}
E\left( X^{n}\right) & =E\left( (P_{1}^{(n)},....,P_{N}^{(n)})\right) \\ 
& =\left( \varphi (P_{1}^{(n)}),....,\varphi (P_{N}^{(n)})\right) \\ 
& =0_{B}=(0,...,0),
\end{array}
$
\end{center}

with

\begin{center}
$P_{j}^{(n)}=\sum_{k=1}^{j}P_{k}^{(n-1)}x_{(j+1)-k},$
\end{center}

where

\begin{center}
$P_{j}^{(1)}=x_{j}$ \ \ and \ $P_{j}^{(2)}=\sum_{k=1}^{j}x_{k}x_{(j+1)-k},$
\end{center}

\strut

in $(A,\varphi ),$ for all $j=1,...,N.$ Hence we have that

\strut

\begin{center}
$\varphi \left( P_{j}^{(n)}\right) =\varphi \left(
\sum_{k=1}^{j}P_{j}^{(n-1)}x_{(j+1)-k}\right) =0\in \Bbb{C},$
\end{center}

\strut

for all $j=1,...,N.$

\strut

(i) Let $n=1.$ Then

\strut

\begin{center}
$\varphi \left( P_{j}^{(1)}\right) =0\in \Bbb{C},$ for $j=1,...,N,$
\end{center}

equivalently,

\begin{center}
$\varphi (x_{j})=0\in \Bbb{C},$ for all $j=1,...,N.$
\end{center}

\strut

(ii) Suppose that $n=2l+1,$ for $l\in \Bbb{N}$. Then

\strut

\begin{center}
$
\begin{array}{lll}
\varphi \left( P_{j}^{(2l+1)}\right) & =\varphi \left(
\sum_{k=1}^{j}P_{j}^{(2l)}x_{(j+1)-k}\right) &  \\ 
& =\sum_{k=1}^{j}\varphi \left( P_{j}^{(2l)}x_{(j+1)-k}\right) =0 & \in \Bbb{%
C},
\end{array}
$
\end{center}

for all $j=1,...,N.$
\end{proof}

\strut

\strut

\begin{example}
Let $N=2.$ Suppose that a $B=\mathcal{C}^{2}$-valued random variable, $%
X=(x_{1},x_{2})\in (T,E),$ is even, where $(T,E)$ is a Toeplitz matricial
probability space over $B,$ induced by a NCPSpace $(A,\varphi ).$ Then we
have the following relation ;

\strut

\begin{center}
$\varphi \left( x_{1}\right) =0=\varphi (x_{2})$
\end{center}

and

\begin{center}
$\varphi \left( P_{1}^{(2l+1)}\right) =\varphi \left( x_{1}^{2l+1}\right) =0$
\ and \ $\sum_{k=1}^{2}\varphi \left( P_{k}^{(2l)}x_{(2+1)-k}\right) =0,$
\end{center}

\strut

in $\Bbb{C},$ for all $l\in \Bbb{N}.$ Note that

\strut

$\sum_{k=1}^{2}\varphi \left( P_{k}^{(2l)}x_{(2+1)-k}\right) =0$

$\ \ \ \ \ \ \ \ \ \ \ \Longleftrightarrow \varphi \left(
P_{1}^{(2l)}x_{2}\right) +\varphi \left( P_{2}^{(2l)}x_{1}\right) =0$

$\ \ \ \ \ \ \ \ \ \ \ \Longleftrightarrow \varphi \left( x_{1}^{2l}\cdot
x_{2}\right) +\varphi \left( P_{2}^{(2l)}x_{1}\right) =0$

$\ \ \ \ \ \ \ \ \ \ \ \Longleftrightarrow \varphi \left(
P_{2}^{(2l)}x_{1}\right) =-\varphi \left( x_{1}^{2l}\cdot x_{2}\right) .$

\strut

Observe that

\begin{center}
$X=(x_{1},x_{2})\equiv \left( 
\begin{array}{ll}
x_{1} & x_{2} \\ 
0 & x_{1}
\end{array}
\right) \in (T,E)$
\end{center}

and

\strut

$X^{2l+1}=(x_{1},x_{2})^{2l+1}=\left( P_{1}^{(2l+1)},\,P_{2}^{(2l+1)}\right) 
$

\strut

$\ \ \ \equiv \left( 
\begin{array}{ll}
P_{1}^{(2l+1)} & P_{2}^{(2l+1)} \\ 
0 & P_{1}^{(2l+1)}
\end{array}
\right) =\left( 
\begin{array}{ll}
x_{1}^{2l+1} & P_{2}^{(2l+1)} \\ 
0 & x_{1}^{2l+1}
\end{array}
\right) =\left( 
\begin{array}{ll}
x_{1}^{2l+1} & P_{2}^{(2l+1)} \\ 
0 & x_{1}^{2l+1}
\end{array}
\right) $

\strut

$\ \ \ =\left( 
\begin{array}{ll}
x_{1}^{2l+1} & x_{1}^{2l}\cdot x_{2}+P_{2}^{(2l)}x_{1} \\ 
0 & \text{ \ \ \ \ \ \ \ \ }x_{1}^{2l+1}
\end{array}
\right) =\left( x_{1}^{2l+1},\,x_{1}^{2l}\cdot x_{2}+P_{2}^{(2l)}\cdot
x_{1}\right) \in (T,E).$

\strut

Therefore, we can get that ;

\strut

(i) \ \ if $a\in (A,\varphi )$ is an (scalar-valued) even random variable
(in the sense that $a\in (A,\varphi )$ is even if and only if all odd
scalar-valued moments vanish), then $\left( a,0\right) \equiv \left( 
\begin{array}{ll}
a & 0 \\ 
0 & a
\end{array}
\right) \in (T,E)$ is again an even $B=\mathcal{C}^{2}$-valued random
variable.

\strut \strut

(ii) \ if $a\in (A,\varphi )$ is an (scalar-valued) even random variable,
then $(a,a)\equiv \left( 
\begin{array}{ll}
a & a \\ 
0 & a
\end{array}
\right) \in (T,E)$ is again an even $B=\mathcal{C}^{2}$-valued random
variable.

\strut

(iii) if $a_{1},a_{2}\in (A,\varphi )$ are \textbf{free} (scalar-valued)
even random variables and if they are identically distributed, then $%
(a_{1},a_{2})\equiv \left( 
\begin{array}{ll}
a_{1} & a_{2} \\ 
0 & a_{1}
\end{array}
\right) \in (T,E)$ is again an even $B=\mathcal{C}^{2}$-valued random
variable, by (ii).
\end{example}

\strut \strut

\begin{corollary}
Let $B=\mathcal{C}^{N}$ be a Toeplitz matricial algebra and let $(A,\varphi
) $ be a NCPSpace in the sense of Section 1.3. Let $(T,E)$ be a Toeplitz
matricial probability space over $B$.

\strut

(i) \ \ If $x\in (A,\varphi )$ is an even (scalar-valued) random variable
(in the sense of Nica, See [2]), then $X=(x,0,...,0)\in (T,E)$ is an even $B$%
-valued random variable.

\strut

(ii) \ If $x\in (A,\varphi )$ is an even (scalar-valued) random variable,
then $X=\left( \underset{N-times}{\underbrace{x,......,x}}\right) \in (T,E)$
is an even $B$-valued random variable.

\strut

(iii) If $x_{1},...,x_{N}\in (A,\varphi )$ are\textbf{\ free} even
(scalar-valued) random variables and if they are identically distributed,
then $X=(x_{1},...,x_{N})\in (T,E)$ is an even $B$-valued random variable.
\end{corollary}

\strut

\begin{proof}
(i) We have that

\strut

\begin{center}
$P_{1}^{(n)}=x^{n}$ and $P_{j}^{(n)}=0_{B}\in (A,\varphi ),$ for all $%
j=2,...,N.$
\end{center}

\strut

Since $x\in (A,\varphi )$ is even, all odd scalar-valued moments vanish. So,

\strut

\begin{center}
$E(X^{n})=E\left( x^{n},0_{B},...,0_{B}\right) $

\strut
\end{center}

$\ \ \ \ \ \ \ \ =\left\{ 
\begin{array}{lll}
\left( \varphi (x^{n}),0,...,0\right) \in B &  & \text{if }n\text{ is even}
\\ 
&  &  \\ 
(0,...,0)=0_{B} &  & \text{if }n\text{ is odd}
\end{array}
\right. $

\strut \strut

$\ \ \ \ \ \ \ \ =\left\{ 
\begin{array}{lll}
E(X^{n}) &  & \text{if }n\text{ is even} \\ 
&  &  \\ 
0_{B} &  & \text{if }n\text{ is odd.}
\end{array}
\right. $

\strut

Thus $X$ is even, in $(T,E).$

\strut

(ii) It is easy to see that

\strut

$X^{n}=\left( x,...,x\right) ^{n}$

\strut

$\ \ \ \ \equiv \left( 
\begin{array}{llll}
x & x & \cdots & x \\ 
& x & \ddots & \vdots \\ 
&  & \ddots & x \\ 
O &  &  & x
\end{array}
\right) ^{n}=\left( 
\begin{array}{llll}
x^{n} & \alpha _{2}x^{n} & \cdots & \alpha _{N}\,x^{n} \\ 
& x^{n} & \ddots & \text{ \ \ }\vdots \\ 
&  & \ddots & \alpha _{2}\,x^{n} \\ 
O &  &  & \text{ \ \ }x^{n}
\end{array}
\right) $

\strut

$\ \ \ \ \equiv \left( x^{n},\alpha _{2}\,x^{n},\alpha
_{3}\,x^{n},...,\alpha _{N}\,x^{n}\right) \in (T,E),$

\strut

where $\alpha _{2},...,\alpha _{N}\in \Bbb{C}$ are arbitrary for each $n\in 
\Bbb{N}.$ So, if $n$ is odd, then

\strut

$E\left( X^{n}\right) =E\left( x^{n},\alpha _{2}\,x^{n},\alpha
_{3}\,x^{n},...,\alpha _{N}\,x^{n}\right) $

\strut

$\ \ \ \ \ \ \ \ \ =\left( \varphi (x^{n}),\,\alpha _{2}\,\cdot \varphi
(x^{n}),...,\alpha _{N}\,\cdot \varphi (x^{n})\right) $

\strut

$\ \ \ \ \ \ \ \ \ =\left\{ 
\begin{array}{ll}
\left( \varphi (x^{n}),\,\alpha _{2}\,\cdot \varphi (x^{n}),...,\alpha
_{N}\,\cdot \varphi (x^{n})\right) & \text{ if }n\text{ is even} \\ 
&  \\ 
\left( 0,...,0\right) =0_{B} & \text{if }n\text{ is odd,}
\end{array}
\right. $

\strut

by the (scalar-valued) evenness of $x\in (A,\varphi ).$ Therefore, $X=\left( 
\underset{N-times}{\underbrace{x,....,x}}\right) \in (T,E)$ is even.

\strut

(iii) Since $x_{1},...,x_{N}$ are free each other, by applying the
M\"{o}bius inversion, we can get the result, similar to (ii).
\end{proof}

\strut

\strut

\strut

\strut

\strut

\strut

\subsection{Toeplitz Matricial Probability Spaces Induced by A Compressed
NCPSpace}

\strut

\strut

\strut

\strut

In [3], Speicher considered a (scalar-valued) compressed NCPSpace and
compressed cumulants. And, in [2], Nica observed the compressed
(scalar-valued) R-transforms. In [32], we observed an opreator-valued
version of them. In this section, we will consider the relation between
Toeplitz matricial probability space $(T_{N},E_{N}),$ induced by a NCPSpace $%
(A,\varphi )$ and Toeplitz matricial probability space $\left(
T_{N}^{p},\,E_{N}^{p}\right) ,$ induced by a compressed NCPSpace $\left( pAp,%
\frac{1}{\varphi (p)}\varphi \mid _{pAp}\right) $ of $(A,\varphi ),$ where $%
p\in A$ is a projection satisfying $\varphi (p)\neq 0.$ By Speicher and
Nica, a compressed NCPSpace $\left( pAp,\varphi _{p}\right) ,$ where $%
\varphi _{p}=\frac{1}{\varphi (p)}\varphi \mid _{pAp},$ is again a NCPSpace,
in the sense of Section 1.3. From now, we will denote a compressed
R-transform of random variable $x\in (A,\varphi )$ and compressed cumulants
of $x,$ as coefficients of the R-transform by $R_{x}^{(pAp)}$ and $%
k_{n}^{(pAp)}(x,...,x),$ respectively. It is known the relation between
cumulants and compressed cumulants, by Speicher and Nica. Also, in [12], if $%
B$ is a unital algebra and $(A,E)$ is a NCPSpace over $B,$ we realized that
there is similar relation between compressed operator-valued cumulants and
operator-valued cumulants, when we take a projection, $p,$ such that $%
E(p)\in C_{A}(B)\cap B_{inv},$ compared with the scalar-valued case.

\strut

\begin{definition}
Let $(A,\varphi )$ be a NCPSPace in the sense of Section 1.3 and let $p\in A$
be a projection (i.e idempotent, $p^{2}=p,$ in $(A,\varphi ),$ in short, $%
p\in A_{pro}$) such that $\varphi (p)\neq 0$. Then $\left( pAp,\varphi
_{p}\right) ,$ with $\varphi _{p}=\frac{1}{\varphi (p)}\varphi \mid _{pAp},$
is called a compressed NCPSpace (by $p\in A_{pro}$). As usual, we will
denote a Toeplitz matricial probability space, induced by $(A,\varphi ),$ by 
$\left( T_{N},E_{N}\right) \equiv (T,E),$ for $N\in \Bbb{N}.$ Denote a
Toeplitz matricial probability space induced by the compressed NCPSpace, $%
\left( pAp,\varphi _{p}\right) ,$ by $\left( T_{N}^{p},\,E_{N}^{p}\right)
\equiv \left( T^{p},E^{p}\right) .$
\end{definition}

\strut

The follwong proposition is proved in [2] and [3] ;

\strut

\begin{proposition}
(See [2] and [3]) Let $(A,\varphi )$ be a NCPSpace and let $p\in A_{pro}$
satisfy $\varphi (p)\neq 0,$ in $\Bbb{C}.$ If $x_{1},...,x_{s}\in (A,\varphi
)$ are (scalar-valued) random variables ($s\in \Bbb{N}$), then

\strut 

\begin{center}
$R_{x_{1},...,x_{s}}^{(pAp)}(z_{1},...,z_{s})=\sum_{n=1}^{\infty }\underset{%
i_{1},...,i_{n}\in \{1,...,s\}}{\sum }\alpha
_{i_{1},...,i_{n}}\,z_{i_{1}}...z_{i_{n}},$
\end{center}

\strut

in $\Theta _{\Bbb{C}}^{s},$ with

\strut

\begin{center}
$
\begin{array}{ll}
\alpha _{i_{1},...,i_{n}} & =k_{n}^{(pAp)}(x_{i_{1}},...,x_{i_{n}}) \\ 
& =\alpha _{0}^{-1}k_{n}\left( \alpha _{0}x_{i_{1}},...,\alpha
_{0}x_{i_{n}}\right) \\ 
& =\alpha _{0}^{n-1}k_{n}(x_{i_{1}},...,x_{i_{n}}),
\end{array}
$
\end{center}

\strut

for all $(i_{1},...,i_{n})\in \{1,...,s\}^{n},$ $n\in \Bbb{N},$ where $%
\alpha _{0}=\varphi (p)\in \Bbb{C}.$ Here, $k_{n}$ is an $n$-th cumulant and 
$k_{n}^{(pAp)}$ is a compressed $n$-th cumulant in the sense of Speicher and
Nica (See [2] and [3]). \ $\square $
\end{proposition}

\strut \strut

\begin{remark}
More generally, in [32], we showed that if $B$ is an arbitrary unital
algebra and $(A,E)$ is a NCPSpace over $B$ and if $p\in A_{pro}$ satisfies
that

\strut

\begin{center}
$E(p)\in C_{A}(B)\cap B_{inv},$
\end{center}

\strut

then $\left( pAp,E_{p}\right) $ is a NCPSpace over $B,$ again, where $%
E_{p}=\left( E(p)\right) ^{-1}E\mid _{pAp}.$ (i.e $B\subset pAp\subset A$)
Now, let $x_{1},...,x_{s}\in (A,E)$ be $B$-valued random variables ($s\in 
\Bbb{N}$) and denote \textbf{trivial} cumulant multiplicative bimodule map
induced by $E$ and those induced by $E_{p}$ by $K_{n}$ and $K_{n}^{(pAp)},$
respectively. Then we have that

\strut

\begin{center}
$
\begin{array}{ll}
K_{n}^{(pAp)}(px_{i_{1}}p,...,px_{i_{n}}p) & =c_{E_{p}}^{(n)}\left(
px_{i_{1}}p\otimes px_{i_{2}}p\otimes ...\otimes px_{i_{n}}p\right) \\ 
& =b_{0}^{-1}c_{E}^{(n)}\left( b_{0}x_{i_{1}}\otimes b_{0}x_{i_{2}}\otimes
...\otimes b_{0}x_{i_{n}}\right) \\ 
& =K_{n}^{symm(b_{0})}\left( x_{i_{1}},...,x_{i_{n}}\right) ,
\end{array}
$
\end{center}

\strut \strut

where $b_{0}=E(p)\in B.$ Therefore, we can get that

\strut

\begin{center}
$R_{px_{1}p,...,px_{s}p}^{(pAp)\,\,:\,%
\,t}(z_{1},...,z_{s})=R_{x_{1},...,x_{s}}^{symm(E(p))}(z_{1},...,z_{s}),$
\end{center}

\strut

in $\Theta _{B}^{s}.$
\end{remark}

\strut

A Toeplitz probability space $(T_{N}^{p},E_{N}^{p})$\strut , induced by a
compressed NCPSpace $(pAp,\varphi _{p}),$ can be understood as an algebra
that

\strut

\begin{center}
$T_{N}^{p}=A\lg \left( \{(px_{1}p,\,px_{2}p,....,px_{N}p)\in
(T_{N},E_{N}):x_{j}\in (A,\varphi ),\text{ }j=1,...,N\}\right) ,$
\end{center}

\strut

where $\left( T_{N},E_{N}\right) $ is a Toeplitz matricial probability space
induced by $(A,\varphi ).$ i.e if $X=(x_{1}^{\prime },...,x_{N}^{\prime
})\in (T_{N}^{p},E_{N}^{p}),$ then there exists $x_{1}^{(i)},...,x_{N}^{(i)}%
\in (A,\varphi ),$ ($i\in \Bbb{N}$) such that

\strut

\begin{center}
$\left( x_{1}^{\prime },...,x_{N}^{\prime }\right) =\sum_{i=1}^{\infty
}\beta _{i}\cdot (px_{1}^{(i)}p,...,px_{N}^{(i)}p),$
\end{center}

\strut

where $\beta _{i}\in \Bbb{C}$ and, for each $i\in \Bbb{N},$

\strut

$\beta _{i}\left( px_{1}^{(i)}p,...,px_{N}p^{(i)}\right) $

$\ \ \ \ \equiv \left( 
\begin{array}{llll}
\beta _{i} &  &  & O \\ 
& \beta _{i} &  &  \\ 
&  & \ddots &  \\ 
O &  &  & \beta _{i}
\end{array}
\right) \left( 
\begin{array}{llll}
px_{1}^{(i)}p & px_{2}^{(i)}p & \cdots & px_{N}^{(i)}p \\ 
& px_{1}^{(i)}p & \ddots & \text{ \ \ }\vdots \\ 
&  & \ddots & px_{2}^{(i)}p \\ 
O &  &  & px_{1}^{(i)}p
\end{array}
\right) $

\strut

$\ \ \ \ =\left( \beta _{i}\cdot px_{1}^{(i)}p,...,\beta _{i}\cdot
px_{N}p^{(i)}\right) .$

\strut

\strut

So, by the above setting, we have that ;

\strut

\begin{theorem}
Let $(T_{N}^{p},E_{N}^{p})$ be a Toeplitz matricial probability space,
induced by a compressed NCPSpace, $\left( pAp,\varphi _{p}\right) ,$ of a
NCPSpace $(A,\varphi ),$ with $p\in A_{pro}$ such that $\varphi (p)\neq 0.$
Let $%
X_{1}=(x_{1}^{(1)},...,x_{N}^{(1)}),...,X_{s}=(x_{1}^{(s)},...,x_{N}^{(s)})%
\in \left( T_{N},E_{N}\right) $ be $\mathcal{C}^{N}$-valued random variables
($s\in \Bbb{N}$) and let $X_{1}^{\prime
}=(px_{1}^{(1)}p,...,px_{N}^{(1)}p),...,X_{s}^{\prime
}=(px_{1}^{(s)}p,...,px_{N}^{(s)}p)\in (T_{N}^{p},E_{N}^{p})$ be $\mathcal{C}%
^{N}$-valued random variables. Then

\strut

$R_{X_{1}^{\prime },...,X_{s}^{\prime }}^{(T_{N}^{p})}(z_{1},...,z_{s})$

\begin{center}
$=\sum_{n=1}^{\infty }\underset{i_{1},...,i_{n}\in \{1,...,s\}}{\sum }\alpha
_{0}^{n-1}K_{n}(X_{i_{1}},...,X_{i_{n}})\,z_{i_{1}}...z_{i_{n}},$
\end{center}

\strut

where $\alpha _{0}=\varphi (p)\in C$ and $K_{n}$'s are cumulant
multiplicative bimodule map induced by $E_{N}:T_{N}\rightarrow \mathcal{C}%
^{N}.$
\end{theorem}

\strut

\begin{proof}
Consider $(i_{1},...,i_{n})$-th coefficient of $R_{X_{1}^{\prime
},...,X_{s}^{\prime }}^{(T_{N}^{p})}$ \ ;

\strut

$coef_{i_{1},...,i_{n}}\left( R_{X_{1}^{\prime },...,X_{s}^{\prime
}}^{(T_{N}^{p})}\right) =K_{n}^{(T_{N}^{p})}\left( X_{i_{1}}^{\prime
},...,X_{i_{n}}^{\prime }\right) $

\strut

$\ \ \ \ =K_{n}^{(T_{N}^{p})}\left(
(px_{1}^{(i_{1})}p,...,px_{N}^{(i_{1})}p),...,(px_{1}^{(i_{n})}p,...,px_{N}^{(i_{n})}p)\right) 
$

\strut

$\ \ \ \ =\left( k_{n}^{(pAp)}\left(
Q_{1\,\,:\,\,p}^{(i_{1},...,i_{n})}\right) ,....,k_{n}^{(pAp)}\left(
Q_{N\,\,\,\,:\,\,\,p}^{(i_{1},...,i_{n})}\right) \right) $

\strut

where $Q_{j\,\,\,\,:\,\,\,p}^{(i_{1},...,i_{n})}\in A^{n}$ is determined
inductively as in Section 2.2 ($j=1,...,N$), with respect to $\left(
T_{N}^{p},\,\,E_{N}^{p}\right) $

\strut

$\ \ \ \ =\left( \alpha
_{0}^{n-1}k_{n}(Q_{1}^{(i_{1},...,i_{n})}),.....,\alpha
_{0}^{n-1}k_{n}(Q_{N}^{(i_{1},...,i_{n})})\right) $

\strut

where $\alpha _{0}=\varphi (p)\in \Bbb{C}$ and $k_{n}$ is a $n$-th cumulant
induced by $\varphi ,$ by the previous proposition. Here, $%
Q_{j}^{(i_{1},...,i_{n})}\in A^{n}$ is determined inductively as in Section
2.2 ($j=1,...,N$), with respect to $(T_{N},\,\,E_{N}),$

\strut

$\ \ \ \ =\alpha _{0}^{n-1}\left(
k_{n}(Q_{1}^{(i_{1},...,i_{n})}),...,k_{n}(Q_{N}^{(i_{1},...,i_{n})})\right) 
$

\strut

by the previous paragraph

\strut

$\ \ \ \ =\alpha _{n}^{n-1}K_{n}\left( X_{i_{1}},...,X_{i_{n}}\right) .$
\end{proof}

\strut

By using notations in [12] and [32], we can get that \ ;

\strut

\begin{corollary}
Under the same hypothesis with the previous theorem,

\strut

\begin{center}
$R_{X_{1}^{\prime },...,X_{s}^{\prime
}}^{(T_{N}^{p})}(z_{1},...,z_{s})=R_{X_{1},...,X_{s}}^{symm(b_{0})}(z_{1},...,z_{s}), 
$
\end{center}

\strut

where

\strut

\begin{center}
$b_{0}=\left( \varphi (p),\,0,....,0\right) \in \mathcal{C}^{N}$
\end{center}

$\square $
\end{corollary}

\strut

\begin{remark}
Let's take a projection $P=(p,0_{A},...,0_{A})\in (T_{N},\,\,E_{N}),$ where $%
p\in (A,\varphi )$ is a projection (i.e $p^{2}=p,$ in $A$). i.e

\strut

\begin{center}
$P=\left( 
\begin{array}{llll}
p &  &  & O \\ 
& p &  &  \\ 
&  & \ddots &  \\ 
O &  &  & p
\end{array}
\right) \in (T_{N},E_{N}),$
\end{center}

\strut

such that $P^{2}=P,$ in $T_{N}.$ Suppose that $\varphi (p)\neq 0.$ Then

\strut

\begin{center}
$E(P)=E\left( p,0_{A},...,0_{A}\right) =\left( \varphi (p),0,...,0\right)
\in \mathcal{C}^{N}.$
\end{center}

\strut

Since $\varphi (p)\neq 0,$ $E(P)\in \mathcal{C}^{N}$ is invertible. So, $%
E(P) $ satisfies the conditions in Remark 3.1. i.e

\strut

\begin{center}
$E(P)\in C_{T_{N}}(\mathcal{C}^{N})\,\,\cap \,\,(\mathcal{C}^{N})_{inv}.$
\end{center}

\strut

And, by the previous corollary, we can see that we have the same result with
Remark 3.1, when we consider the Toeplitz matricial probability space $%
(T_{N}^{p},\,\,E_{N}^{p}),$ induced by $(pAp,\varphi _{p}).$ i.e this
Toeplitz matricial probability space $(T_{N}^{p},\,\,E_{N}^{p})$ is nothing
but a compressed (amalgamated) NCPSpace over $\mathcal{C}^{N},$ $\left(
PT_{N}P,\,\,E_{N}(P)^{-1}E_{N}\mid _{PT_{N}P}\right) .$
\end{remark}

\strut \strut

\begin{corollary}
Let $(A,\varphi )$ be a NCPSpace, in the sense of Section 1.3 and let $%
\left( pAp,\varphi _{p}\right) $ be a compressed NCPSpace, by $p\in A_{pro},$
with $\varphi (p)\neq 0,$ in $\Bbb{C}.$ Let $(T_{N},\,\,E_{N})$ and $\left(
T_{N}^{p},\,\,E_{N}^{p}\right) $ be Toeplitz matricial probability spaces
induced by $(A,\varphi )$ and $\left( pAp,\varphi _{p}\right) ,$
respectively. Then $\left( T_{N}^{p},\,\,E_{N}^{p}\right) $ is an
amalgamated compressed NCPSpace over $\mathcal{C}^{N},$ by $P=\left(
p,0_{A},...,0_{A}\right) \in (T_{N},\,\,E_{N}).$ i.e

\strut

\begin{center}
$\left( T_{N}^{p},\,\,E_{N}^{p}\right) =\left(
PT_{N}P,\,\,E_{N}(P)^{-1}E_{N}\mid _{PT_{N}P}\right) .$
\end{center}

$\square $
\end{corollary}

\strut

\strut

A Toeplitz matricial probability theory can be a good example of amalgamated
free probability theory (See [1] and [12]). Futhermore, since $C_{T_{N}}(%
\mathcal{C}^{N})=$ \strut $\mathcal{C}^{N},$ it is easy to deal with. i.e we
can avoid a difficulty caused by Insertion property in Amalgamated free
probability and hence it is enough to consider trivial R-transforms, in the
sense of Section 1.2, (defined as Toeplitz matricial R-transforms, in this
paper) when we work for R-transform theory and R-transform calculus, on
Toeplitz matricial probability spaces.

\strut

\strut \strut

\strut

\strut

\strut

\strut

\strut

\strut

\strut

\strut

\strut

\strut \textbf{References}

\strut

\label{REF}

\strut

\begin{enumerate}
\item[{[1]}]  R. Speicher, Combinatorial Theory of the Free Product with
Amalgamation and Operator-Valued Free Probability Theory, AMS Mem, Vol 132 ,
Num 627 , (1998).

\strut

\item[{[2]}]  A. Nica, R-transform in Free Probability, IHP course note,
available at

\texttt{www.math.uwaterloo.ca/\symbol{126}anica.}

\strut

\item[{[3]}]  R. Speicher, Combinatorics of Free Probability Theory IHP
course note,

available at \texttt{www.mast.queensu.ca/\symbol{126}speicher.}

\strut

\item[{[4]}]  A. Nica, D. Shlyakhtenko and R. Speicher, R-cyclic Families of
Matrices in Free Probability, J. of Funct Anal, 188 (2002), 227-271.

\strut

\item[{[5]}]  A. Nica and R. Speicher, R-diagonal Pair-A Common Approach to
Haar Unitaries and Circular Elements, (1995), \texttt{www.mast.queensu.ca/%
\symbol{126}speicher.}

\strut

\item[{[6]}]  D. Shlyakhtenko, Some Applications of Freeness with
Amalgamation, J. Reine Angew. Math, 500 (1998), 191-212.

\strut

\item[{[7]}]  A. Nica, D. Shlyakhtenko and R. Speicher, R-diagonal Elements
and Freeness with Amalgamation, Canad. J. Math. Vol 53, Num 2, (2001)
355-381.

\strut

\item[{[8]}]  A. Nica, R-transforms of Free Joint Distributions and
Non-crossing Partitions, J. of Func. Anal, 135 (1996), 271-296.

\strut

\item[{[9]}]  D.Voiculescu, K. Dykemma and A. Nica, Free Random Variables,
CRM Monograph Series Vol 1 (1992).

\strut

\item[{[10]}]  D. Voiculescu, Operations on Certain Non-commuting
Operator-Valued Random Variables, Ast\'{e}risque, 232 (1995), 243-275.

\strut

\item[{[11]}]  D. Shlyakhtenko, A-Valued Semicircular Systems, J. of Funct
Anal, 166 (1999), 1-47.

\strut

\item[{[12]}]  I. Cho, Operator-Valued Free Probability and its Application
(preprint).

\strut

\item[{[13]}]  R. Speicher, A Conceptual Proof of a Basic Result in the
Combinatorial Approach to Freeness, \texttt{www.mast.queensu.ca/\symbol{126}%
speicher.}

\strut

\item[{[14]}]  \strut M. Bo\.{z}ejko, M. Leinert and R. Speicher,
Convolution and Limit Theorems for Conditionally Free Random Variables, 
\texttt{www.mast.queensu.ca/\symbol{126}speicher.}
\end{enumerate}

\strut

\begin{enumerate}
\item[{\lbrack 15]}]  A. Nica, R-diagonal Pairs Arising as Free Off-diagonal
Compressions, available at

\texttt{www.math.uwaterloo.ca/\symbol{126}anica.}
\end{enumerate}

\strut \strut

\begin{enumerate}
\item[{\lbrack 16]}]  B.Krawczyk and R.Speicher, Combinatorics of Free
Cumulants, available at

\texttt{www.mast.queensu.ca/\symbol{126}speicher.}
\end{enumerate}

\strut \strut

\begin{enumerate}
\item[{\lbrack 17]}]  A. Nica and R.Speicher, A ''Fouries Transform'' for
Multiplicative Functions on Noncrossing Patitions, J. of Algebraic
Combinatorics, 6, (1997) 141-160.
\end{enumerate}

\strut

\begin{enumerate}
\item[{\lbrack 18]}]  P.\'{S}niady and R.Speicher, Continous Family of
Invariant Subspaces for R-diagonal Operators, Invent Math, 146, (2001)
329-363.
\end{enumerate}

\strut

\begin{enumerate}
\item[{\lbrack 19]}]  I.Cho, Amalgamated S-transforms, preprint.
\end{enumerate}

\strut

\begin{enumerate}
\item[{\lbrack 20]}]  I.Cho, $B$-Poisson distribution and Amalgamated
Conditional Freeness, preprint.
\end{enumerate}

\strut

\begin{enumerate}
\item[{\lbrack 21]}]  Y. Watatani, Index for C$^{*}$-subalgebras, Mem of
AMS, Vol 83, 424 (1990).
\end{enumerate}

\strut

\begin{enumerate}
\item[{\lbrack 22]}]  V. Jones, Index for Subfactors, Invent. Math, 72
(1983) 1-25.
\end{enumerate}

\strut

\begin{enumerate}
\item[{\lbrack 23]}]  V. Jones, Braid Groups, Hecke Algebras and Type II$%
_{1} $-factors, Geometric Methods in Operator Algebra, Pitman Res. Notes 123
(1986) 242-273.
\end{enumerate}

\strut

\begin{enumerate}
\item[{\lbrack 24]}]  V. Jones, Hecke Algebra Representations of Braid
Groups and Link Polynomials, Ann Math, 126 (1987) 335-388.
\end{enumerate}

\strut

\begin{enumerate}
\item[{\lbrack 25]}]  H. Kosaki, Extension of Jones' Theory on Index to
Arbitrary Factors, J. of F.A, 66 (1986), 123-140.
\end{enumerate}

\strut

\begin{enumerate}
\item[{\lbrack 26]}]  M. Pimsner and S. Popa, Iterating the Basic
Construction, Trans. AMS, 310, No 1 (1988), 127-133.
\end{enumerate}

\strut

\begin{enumerate}
\item[{\lbrack 26]}]  S. Popa, Markov Traces on Universal Jones Algebras and
Subfactors of Finite Index. Invent. Math, 111 (1993), 375-405.
\end{enumerate}

\strut

\begin{enumerate}
\item[{\lbrack 27]}]  I. Cho, Amalgamated R-transform Theory of Index-Finite
type Conditional Expectations, preprint.
\end{enumerate}

\strut

\begin{enumerate}
\item[{\lbrack 28]}]  I. Cho, $B$-free family of Cuntz Isometries, preprint.
\end{enumerate}

\strut

\begin{enumerate}
\item[{\lbrack 29]}]  I. Cho, Amalagmated Free Probability Theory and its
Application, preprint.
\end{enumerate}

\strut

\begin{enumerate}
\item[{\lbrack 30]}]  A. Nica, D. Shlyakhtenko, R. Speicher, Operator-Valued
Distributions I. Characterizations of Freeness, preprint.
\end{enumerate}

\strut

\begin{enumerate}
\item[{\lbrack 31]}]  A. Nica, D. Shlyakhtenko, F. Goodman, Free Probability
of Type B, preprint.
\end{enumerate}

\strut

\begin{enumerate}
\item[{\lbrack 32]}]  I. Cho, Compressed Amalgamated R-transform Theory,
preprint.
\end{enumerate}

\strut

\begin{enumerate}
\item[{\lbrack 33]}]  I. Cho, Perturbed R-transform Theory, preprint.
\end{enumerate}

\strut

\end{document}